\documentclass[10pt,a4paper,reqno]{amsart}

\usepackage[english]{babel}
\usepackage{amsfonts, amsmath, amsthm, amssymb,amscd,indentfirst}

\usepackage{graphicx}
\usepackage{times}
\usepackage{palatino}
\usepackage{tikz}
\usepackage{marginnote}
\usepackage{stackrel}
\usepackage{lipsum}

\usepackage{float}
\usepackage{enumerate}

\usepackage{anysize}
\marginsize{2.9cm}{2.9cm}{2.5cm}{2.5cm}

\numberwithin{equation}{section}

\newtheorem{theorem}{Theorem}[section]
\newtheorem{proposition}[theorem]{Proposition}
\newtheorem{lemma}[theorem]{Lemma}
\newtheorem{definition}[theorem]{Definition}
\newtheorem{corollary}[theorem]{Corollary}

\newtheorem{problem}[theorem]{Problem}

\newtheorem{ex}[theorem]{Example}
\newtheorem{remark}[theorem]{Remark}

\newtheorem{assumption}{Assumption}[section]

\newtheorem{cl}{Claim}
\newtheorem*{cl*}{Claim}

\newcommand{\dv}{\textnormal{d}v}
\newcommand{\s}{\textnormal{d}\sigma}
\newcommand{\dx}{\textnormal{d}x}

\title{A singular Yamabe problem on manifolds with solid cones}
\author{Juan Alcon Apaza}
\author{Sergio Almaraz}

\address{Juan Alcon Apaza, Universidade Federal Fluminense, Instituto de Matemática, Campus do Gragoatá Rua Prof. Marcos Waldemar de Freitas, s/n, bloco H, Niterói, RJ 24210-201, Brazil}
\email{jpablo@id.uff.br}

\address{Sergio Almaraz, Universidade Federal Fluminense, Instituto de Matemática, Campus do Gragoatá Rua Prof. Marcos Waldemar de Freitas, s/n, bloco H, Niterói, RJ 24210-201, Brazil}
\email{sergioalmaraz@id.uff.br}

\begin{document}

\maketitle

\begin{abstract}
We study the existence of conformal metrics on noncompact Riemannian manifolds with noncompact boundary, 
which are complete as metric spaces and have negative constant scalar curvature in the interior and negative constant mean curvature on the boundary. 
These metrics are constructed on smooth manifolds obtained by removing $d$-dimensional submanifolds from certain n-dimensional 
compact spaces locally modelled on generalized solid cones. We prove the existence of such metrics if and only if $d>(n-2)/2$. 
Our main theorem is inspired by the classical results by Aviles-McOwen and Loewner-Nirenberg known in the literature as the $``$singular Yamabe problem$”$. 
\end{abstract}



\let\thefootnote\relax
\let\thefootnote\relax
\let\thefootnote\relax\footnote{The first author was supported by CAPES-88882.456643/2019-01. The second author was partially supported by FAPERJ-202.802/2019.}

\section{Introduction}\label{intro}

The singular Yamabe problem is an extension, to noncompact Riemannian manifolds, of the classical Yamabe problem.
While the classical problem \cite{Y} is formulated for closed manifolds (i.e., compact and without boundary), its singular version concerns manifolds obtained by removing a closed subset from a closed manifold. Precisely, the following problem is considered:
\begin{problem}\label{problema1} 
Let  $(M,g)$ be a closed smooth Riemannian manifold of dimension $n\geq 3$ and let $F\subset M$ be a closed subset (which we call a $``$singular set$"$). Is there a complete metric on $M\backslash F$ which is conformal to $g$ and has scalar curvature $R=\operatorname{constant}$?
\end{problem} 

As the following example suggest, the dimension of the singular set is closely related to the sign of $R$:

\begin{ex}\label{canonico:ex}
Let $\mathbb{S}^n$ be the $n$-dimensional unit sphere and let $\mathbb{S}^d$ be a  sphere  of  dimension $0\leq d\leq n-1$, totally geodesic in $\mathbb{S}^n$. The stereographic projection $\psi:\mathbb{S}^n\backslash \{p\}\to \mathbb R^n$ gives a conformal equivalence between $\mathbb{S}^n\backslash \mathbb{S}^d$ and  $\mathbb R^n\backslash\mathbb R^d$, where $p$ is the north pole of $\mathbb{S}^n$ and
$$
\mathbb R^d=\{x=(x_1,...,x_n)\in\mathbb R^n\:|\:x_{d+1}=...=x_n=0\}.
$$
On the other hand, the conformal change 
$$
dx_1^2+...+dx_n^2\mapsto \frac{dx_1^2+...+dx_n^2}{x_{d+1}^2+...+x_n^2}
$$ 
transforms  $\mathbb R^n\backslash\mathbb R^d$ into $\mathbb H^{d+1}\times  \mathbb{S}^{n-d-1}$, where $\mathbb H^{k}$ represents the hyperbolic space of dimension $k$. This product manifold has constant scalar curvature
$$
R=(n-d-1)(n-d-2)-d(d+1)=(n-1)(n-2-2d),
$$
whose signal is determined by $\frac{n-2}{2}-d$. 
\end{ex}

Historically, the Singular Yamabe problem originated from the seminal paper \cite{LN}, by Loewner and Nirenberg, which handled the case $R<0$ on spheres. 
The case of a general compact Riemannian manifolds was studied by Aviles and McOwen:

\begin{theorem} {\cite{aviles}} \label{148}
Suppose that $(M,g)$ is a closed Riemannian manifold of dimension $n$. Let $F$ be a closed smooth submanifold of $M$, with dimension $0\leq d\leq n-1$. Then the  Problem \ref{problema1} has a solution with  $R=\operatorname{constant}<0$ if and only if $d>\frac{n-2}{2}$. 
\end{theorem}

An immediate consequence of the particular case $d=n-1$ is stated as follows:

\begin{corollary}\label{corol:bdry}
If $(M,g)$ is a compact Riemannian manifold with non-empty  smooth boundary $\partial M$, then $M\backslash\partial M$ admits a complete metric, conformal to $g$, with constant negative scalar curvature.
\end{corollary}

\begin{remark}
We refer the reader to \cite{finn1, finn2, labutin} for later developments considering more general singular sets. 
The corresponding fully nonlinear version of the singular Yamabe problem was studied in \cite{gonz}.
\end{remark}

It is natural to seek for an extension, to noncompact manifolds, of the Escobar-Yamabe problem for compact manifolds with boundary introduced in \cite{escobar}, where one searches for a conformal metric with constant scalar curvature in the interior and constant mean curvature on the boundary. 
The first challenge we face is the quest for a compact space from which one should remove a singular subset. It turns out that the intersection angle between the boundary and the singular set plays an important role in the analysis as the next example suggests.

\begin{ex}\label{model:ex}
For $h>0$ and $d\in\{1,...,n-1\}$, consider the generalized solid cone 
$$
\mathcal{C} _{d,h} = \left\{x\in\mathbb R^n\:|\:x_1\geq h \left(x_{d+1} ^2 + \cdots + x_n ^2 \right) ^{\frac{1}{2}}\right\},
$$
and the singular set 
$$
\mathbb R^d_+=\{x\in\mathbb R^n\:|\:x_1\geq 0,\,x_{d+1}=...=x_n=0\},
$$
as in Figure \ref{solidcone}.
\begin{figure}[H]
\centering
\includegraphics[scale=0.17]{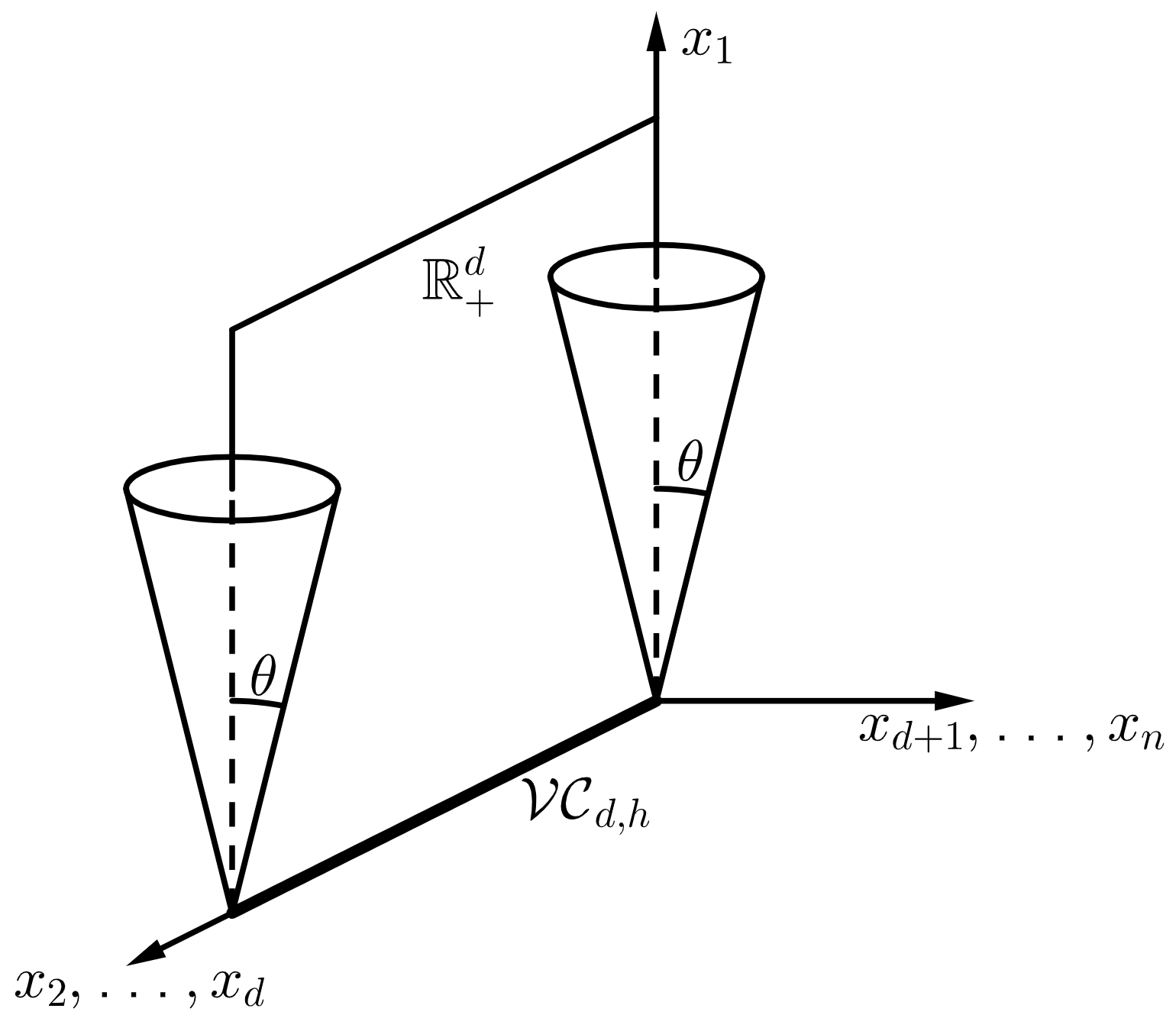}
\caption{The generalized solid cone $\mathcal{C} _{d,h}$.}
\label{solidcone}
\end{figure}

Observe that $\mathbb R^d_+$ intersects the boundary cone
$$
\partial \mathcal{C} _{d,h} = \left\{x\in\mathbb R^n\:|\:x_1= h \left(x_{d+1} ^2 + \cdots + x_n ^2 \right) ^{\frac{1}{2}}\right\},
$$
along the $``$vertex$"$ 
$$
\mathcal V\mathcal{C} _{d,h}=\{x\in\mathbb R^n\:|\:x_1=0=x_{d+1}=...=x_n\},
$$
with an angle $\theta\in(0,\pi/2)$ determined by $h=\cot\theta$.
Now, $\mathcal{C} _{d,h}\backslash\mathbb R^d_+\subset \mathbb R^n\backslash \mathbb R^d$ is a manifold with smooth noncompact boundary 
$$
\partial (\mathcal{C} _{d,h}\backslash\mathbb R^d_+)= \left\{ x\in\mathbb R^n\:|\:x_1= h \left(x_{d+1} ^2 + \cdots + x_n ^2 \right) ^{\frac{1}{2}}, \ x_1>0 \right\}
$$
in the classical sense. Hence, $\mathcal{C}_{d,h} \backslash \mathbb{R}^d _+$, with the metric
\begin{equation}\label{model:metric}
\frac{dx_1^2+...+dx_n^2}{x_{d+1}^2+...+x_n^2}
\end{equation}
is a region of $\mathbb H^{d+1}\times \mathbb S^{n-d-1}$, so that it has scalar curvature $(n-1)(n-2-2d)$.

On the other hand, a direct calculation carried out in Section \ref{155}   shows that $\partial (\mathcal{C} _{d,h}\backslash\mathbb R^d_+)$ has constant mean curvature 
$$
\frac{-dh}{\sqrt{h^2+1}}=-d\cos\theta.
$$
Observe that, while the sign of the scalar curvature is determined by $\frac{n-2}{2}-d$, the sign of the boundary mean curvature is determined by $-\cos\theta$.

Adding the infinity of $\mathbb R^n$ to that structure, one obtains the compact topological space
$$
M=\mathcal C_{d,h}\cup\{\infty\}.
$$
If we set $\Gamma =\mathbb R^d_+\cup\{\infty\}$ we obtain the noncompact manifold $M\backslash\Gamma=\mathcal C_{d,h}\backslash \mathbb R^d_+$ with noncompact boundary $\partial(\mathcal C_{d,h}\backslash\mathbb R^d_+)$ which will be our model in this paper. 
So, in analogy with Problem \ref{problema1}, the pair $(M,\Gamma)$ plays the same role as $(\mathbb S^n, \mathbb S^d)$ does in Example \ref{canonico:ex}. The reader may find interesting to see $M$ as a compact subset of $\mathbb S^n$ by means of the stereographic projection; see Figure \ref{sphere}.

\begin{figure}[H]
\centering
\includegraphics[scale=0.15]{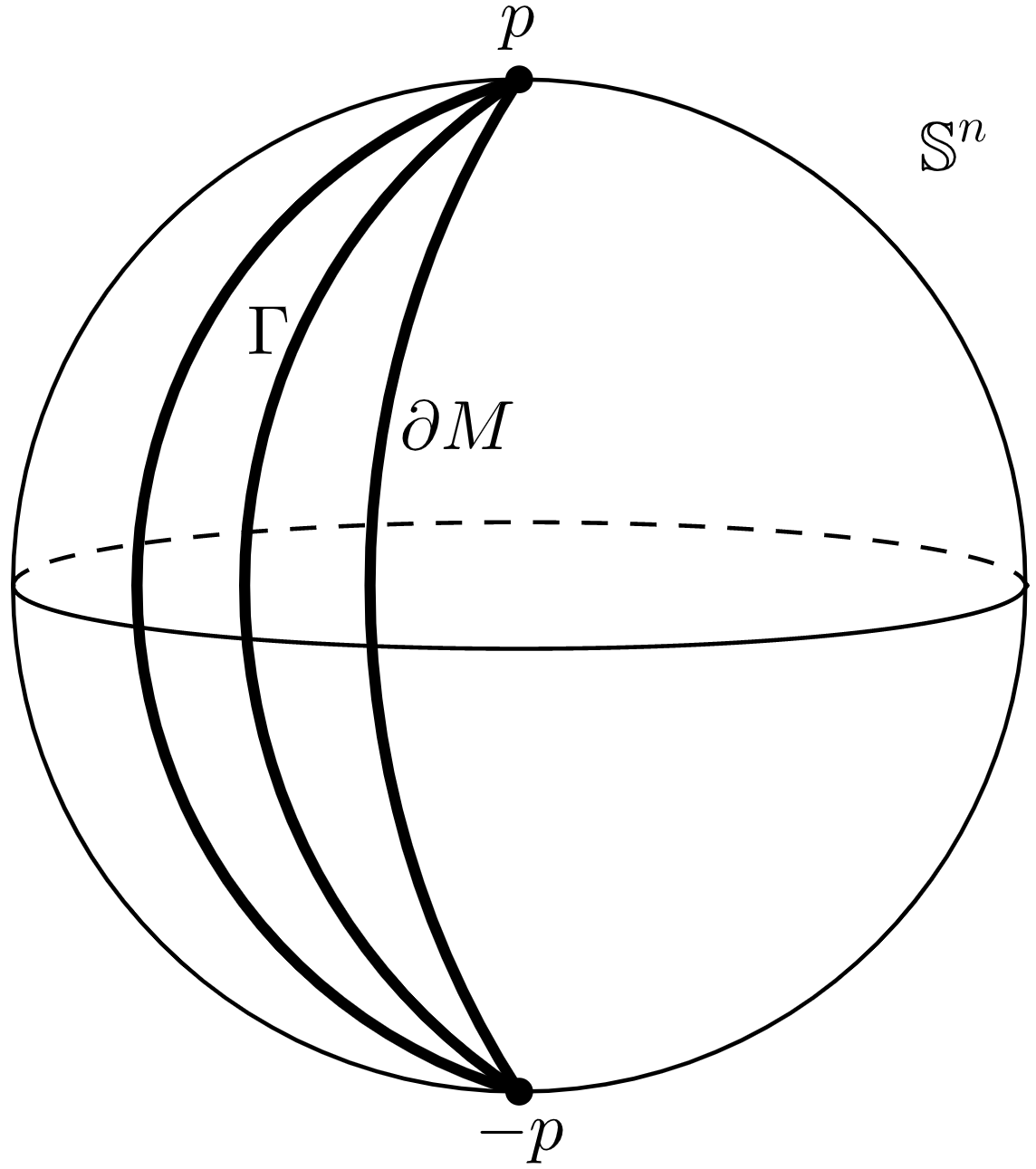}
\caption{The model example under the stereographic projection.}
\label{sphere}
\end{figure}
\end{ex}

The above construction motivates us to define a $``$solid conical manifold$"$ as a topological space $M$ with a differential structure locally modeled on $\mathcal C_{d,h}$; see Section \ref{155}   for the precise definitions. This produces a nonsmooth boundary $\partial M$ corresponding to the points of $M$ taken to $\partial \mathcal C_{d,h}$ by the coordinate charts, and a vertex $\mathcal V M$ corresponding to the ones taken to $\mathcal V\mathcal C_{d,h}$.
From one such solid conical manifold, we remove a closed singular set $\Gamma$ which is given by a finite union $\cup_m\Gamma_m$  of $d_m$-dimensional submanifolds $\Gamma_m$ with boundary $\partial \Gamma_m$, in such a way that $\mathcal V M=\cup_m\partial\Gamma_m$.

To that structure we add a smooth Riemannian metric $g$ on $M$. This metric defines the intersection angle between $\Gamma$ and $\partial M$ along the vertex $\mathcal{V}M$, which is assumed to coincide with the angle $\theta$ coming from the solid cone modeling it; see Definition \ref{107a}.


We raise the following question: 
\begin{problem}\label{problem2}
Is there a Riemannian metric on $M\backslash \Gamma$, complete as a metric space and conformal to $g$, that has  constant scalar curvature  in $M\backslash \Gamma$  and  constant  mean curvature on $\partial (M\backslash \Gamma)=\partial M \backslash \mathcal V M$?
\end{problem}

One important feature of this problem is that, although $M$ is not a manifold in the classical sense, $M\backslash \Gamma$ is a noncompact manifold with smooth noncompact boundary $\partial(M\backslash \Gamma)=\partial M\backslash \mathcal V M$ in the classical sense. So, this can be viewed as an existence problem of conformal metrics on noncompact Riemannian manifolds with constant scalar curvature and constant boundary mean curvature, while $M$ is viewed as the compactification of $M\backslash \Gamma$.

Our approach to Problem \ref{problem2} is a modification of the arguments in \cite{aviles, 1aviles} which depend essentially on maximum principles, variational technics and elliptical estimates adapted to our settings.
In analytical terms, we are searching for a solution $u>0$ of the problem
\begin{equation*}\label{sing:escobar}
\left\{
\begin{aligned}
-\Delta _g u+\displaystyle\frac{(n-2)}{4(n-1)} R_g u+c_0 u^{\frac{n+2}{n-2}}&=0   & & \text { in } M\backslash \Gamma ,\\
\displaystyle\frac{\partial u}{\partial \nu _g}+\displaystyle\frac{n-2}{2(n-1)} H_g u+c_1 u^{\frac{n}{n-2}}&=0 & & \text { on } \partial M \backslash \mathcal V M ,\\
\liminf_{p \rightarrow \Gamma} u(p)\operatorname{dist}_g(p, \Gamma)^{\frac{n-2}{2}} &>0.
\end{aligned}
\right.
\end{equation*}
Here, $c_0>0$ and $c_1>0$ are constant on each connected component of $M\backslash \Gamma$ and $\partial M \backslash \mathcal V M$  respectively. Besides, $R_g$ stands for the scalar curvature, $\Delta_g=\operatorname{div}_g\nabla_g$ is the Laplace operator and  $\nu _g$ is the outward unit normal vector to $\partial M\backslash  \mathcal V M$, so that $H_g=\operatorname{div}_g\nu _g$ is its mean curvature.

Our main result, see Theorem \ref{209} below, implies in particular that, given $c_0, c_1>0$ as above, one can find a solution $\widetilde g$ to the Problem \ref{problem2} with
$$
\left\{
\begin{aligned}
R_{\tilde g}&=-c_0, & &  M\backslash\Gamma ,\\
H_{\tilde g}&=-c_1, & & \partial M \backslash \mathcal V M,
\end{aligned}
\right.
$$
if and only if $\operatorname{dim} \Gamma_m>\frac{n-2}{2}$ for all $m$.
In fact, our result is more general as it allows for $c_0$ and $c_1$ to be smooth functions on $ M \backslash \Gamma$ and $\partial M \backslash \mathcal V M$ respectively, bounded above and below by positive constants.

\begin{figure}[H]
\centering
\includegraphics[scale=0.14]{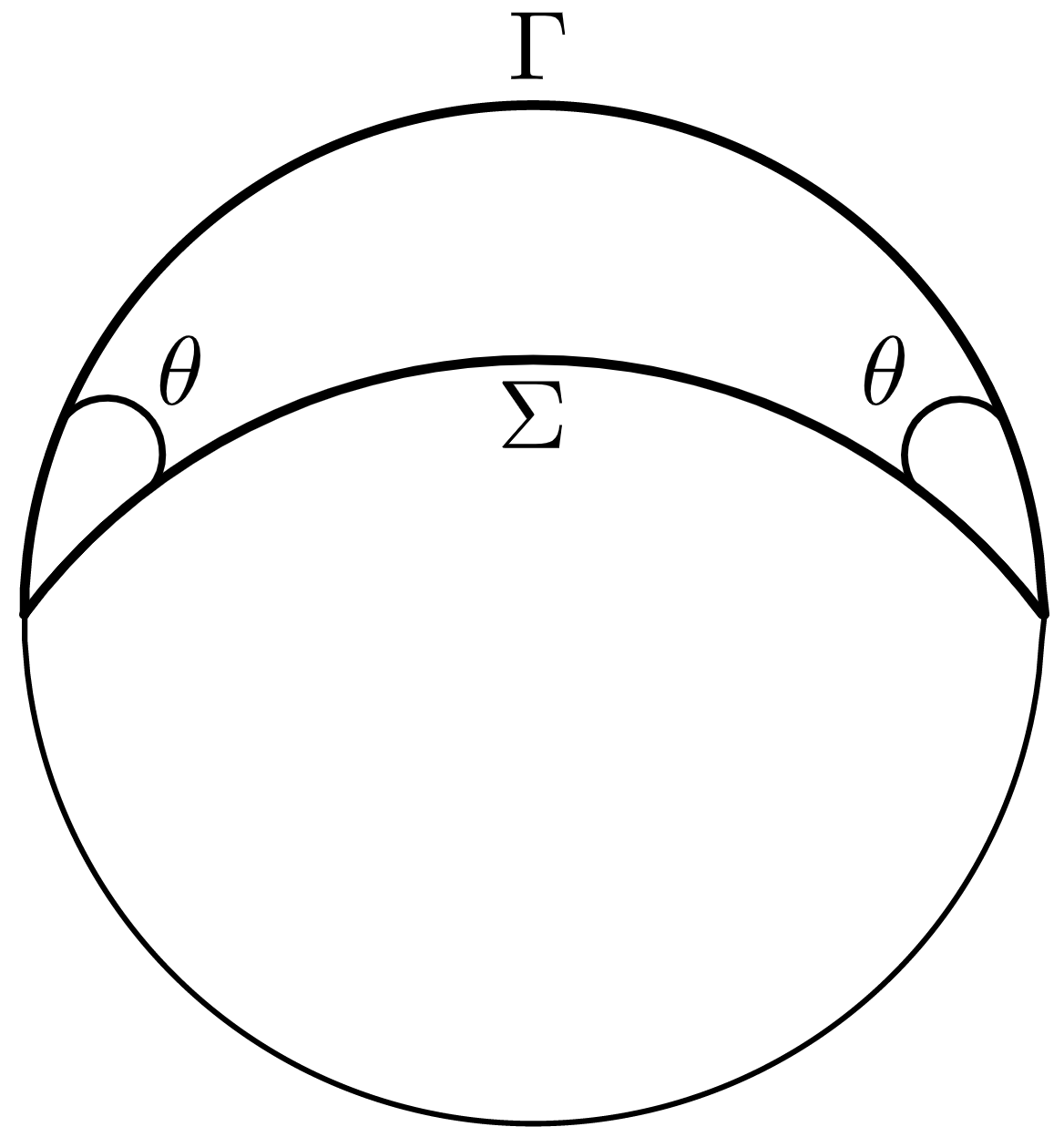}
\caption{One connected  component  of $\mathcal{C} _{h,n-1} \cup \{\infty\}$ in the Poincaré ball model for the hyperbolic space (recall that $h=\operatorname{cot} \theta $).}
\label{d9}
\end{figure}

\begin{remark}\label{rmk:corner}
The case $d=n-1$ has special interest because of its close relationship with cornered manifolds (see Appendix \ref{105}). In our model, Example \ref{model:ex}, the singular set 
$$
\mathbb R^{n-1}_+=\{x\in\mathbb R^n\:|\:x_1\geq 0=x_n\},
$$
splits $M=\mathcal C_{n-1,h}\cup\{\infty\}$ in two connected components, each one being  a cornered manifold itself. In this case, the metric (\ref{model:metric}) is the hyperbolic metric 
$$
x_n^{-2}(dx_1^2+...+dx_n^2),
$$
so that each connected component of $\mathbb R^n\backslash \mathbb R^{n-1}$ is a copy of $\mathbb H^n$. Further,  the two connected components of the boundary 
$x_1=h|x_n|\neq 0$  are hyperspheres of $\mathbb H^n$, i.e., hypersurfaces equidistant from the totally geodesic hypersurface $x_1 = 0$; see Figure \ref{d9}. 

More generally, let $M$ be a cornered manifold $M$ with boundary $\partial M=\Gamma\cup\Sigma$, where $\Gamma$ and $\Sigma$ are smooth hypersurfaces of $M$ with the common boundary $\partial\Gamma=\partial\Sigma=\Gamma\cap\Sigma$ being a codimension two corner of $M$. 
Let $g$ be a Riemannian metric on $M$. If we assume that  $\Gamma$ and $\Sigma$ make intersection angles $\theta_m$ along each connected component of $\Gamma\cap\Sigma$, then doubling $M$ along $\Gamma$ produces a solid conical manifold with $d_m=n-1$ and $h_m=\cot \theta_m$. 
\end{remark}

The same proof of Theorem \ref{209} also leads to the following analogue of Corollary \ref{corol:bdry}:
\begin{corollary}\label{corol:corner}
Let $M$ be a cornered manifold as in Remark \ref{rmk:corner} and let $g$ be a Riemannian metric on $M$. Assume that the angle between $\Gamma$ and $\Sigma$ along each connected component of the corner $\Gamma\cap\Sigma$ is a constant in the interval $(0,\pi/2)$.
Given smooth functions $c_0$ and $c_1$ on $M\backslash \Gamma$ and $\Sigma$ respectively, both being bounded above and below by positive constants, there exists a smooth metric $\widetilde g\in [g]$ on $M\backslash \Gamma$, complete as a metric space, satisfying
$$
\left\{
\begin{aligned}
R_{\tilde g}&=-c_0, &  & M\backslash \Gamma ,\\
H_{\tilde g}&=-c_1, & & \Sigma.
\end{aligned}
\right.
$$
\end{corollary}

We conclude this section by discussing a couple of questions arising naturally from our results.  
The first one is the uniqueness of the metrics obtained in Corollary \ref{corol:corner}. In the case of Corollary \ref{corol:bdry}, uniqueness was proved in \cite{ander} and \cite{maz}, where regularity of solutions up to the singular set is studied. 
We believe that similar results should hold in our settings, where an asymptotic behavior is expected with respect to the boundary component $\Gamma$. We also believe that a hyperboloidal initial data problem for Einstein's field equations could be formulated in the presence of a noncompact boundary. 
This motivates one to pursue results similar to \cite{ander}.
In that case, the boundary mean curvature would be coupled with the scalar in the interior to compose the hyperboloidal initial data set. This is strongly motivated by the results in \cite{AdLM} where the boundary mean curvature plays central role in the definitions of dominant energy conditions and total mass invariants.

The second question is the case when $c_1$ and $c_2$ are not both negative. Although the case $c_1=0$ and $c_2 <0$ should be similar to the one when $c_1, c_2<0$, handled here, the situation when $c_1$ is positive should require more refined techniques as the behavior near the singular set is not a priori determined. 
We refer the reader to \cite{byde, mazzeo-pacard1, mazzeo-pacard2, mazzeo-pollack-uhlenbeck, mazzeo-smale, santos, schoen} for interesting results that could be extended to our setting.

This paper is organized as follows.  In Section \ref{155}, we give precise definitions of the objects involved and discuss further our model, namely, Example  \ref{model:ex}. In  Section \ref{210}, we state and prove our main  result, Theorem \ref{209}. The Appendix contains some technical tools used in this paper.


\section{Preliminaries and formal definitions} \label{155}

In this section, we define the type of manifolds appearing in our main theorem stated and proved in Section \ref{210} below. For $n\geq 2$, they are modeled on sets of the form
$$
\mathcal{C} _{d,h}=\left\{\left(x_{1}, \ldots, x_{n}\right) \in \mathbb{R} ^{n} \:|\: x_{1} \geq h \left(x_{d+1}^2 + \cdots + x_{n}^2\right)^{\frac{1}{2}}\right\},
$$
where $h> 0$ and $d\in \{1, \ldots, n-1\}$, as described in the Introduction. We now provide the precise definitions.

\begin{definition} A \textnormal{smooth solid conical manifold} of dimension  $n$  is a paracompact Hausdorff topological space  $M$  and a family of homeomorphisms    (called the charts)
$$
\varphi_{\alpha}: U_{\alpha} \rightarrow \varphi_{\alpha}(U_{\alpha})\subset M,
$$ 
where each $U_{\alpha} \subset \mathcal C _{d_{\alpha} , h_{\alpha}}$ is a relative open subset, satisfying the following three conditions:
\begin{enumerate}[$(i)$]
\item $\bigcup_{\alpha} \varphi_{\alpha}\left( U _{\alpha}\right)=M$.
\item For any pair  $\alpha$, $\beta$,   with  $W:= \varphi_{\alpha}\left( U _{\alpha}\right) \cap \varphi_{\beta}\left( U _{\beta}\right) \neq \varnothing$,  the map  
$$
\varphi_{\beta}^{-1} \circ \varphi_{\alpha}: \varphi_{\alpha}^{-1}(W) \rightarrow \varphi_{\beta}^{-1}(W)
$$  
is smooth ${}^1 $\footnote{${}^1$ Suppose  $A \subset \mathbb{R}^{n}$  and  $B\subset\mathbb{R}^m$,   we say that a map  $f: A \rightarrow B$  is smooth if, for any $p\in A$, $f$ has a smooth extension in a neighborhood of  $p$, i.e., there exist an open subset  $V \subset \mathbb{R} ^n$  with  $p\in V$  and  a smooth map  $ \widetilde{f}: V \rightarrow \mathbb{R} ^m$  with   $\widetilde{f} |_{V\cap A} \equiv f|_{V\cap A}$.}. 
\item The family $\left\{\left(U_{\alpha}, \varphi_{\alpha}\right)\right\}$ is maximal relative to the conditions $(i)$  and $(ii)$.
\end{enumerate}
\end{definition}

We define the \textit{vertex} 
$$
\mathcal{V} M=\left\{p \in M \:|\:  \exists \left(U_\alpha, \varphi_\alpha\right)  \text { and }  x \in U_\alpha  \text { such that }  x_{1}= 0 \text { and }  \varphi _{\alpha}(x)=p \right\},
$$
the \textit{interior}
$$
\operatorname{int} (M)=\left\{p \in M \:|\: \exists \left(U_\alpha, \varphi_\alpha\right)   \text { such that }   p\in \varphi _{\alpha} \left( U_{\alpha}\right) 
 \text { and } U_\alpha \text { is an open subset of } \mathbb{R}^n \right\},
$$
and the \textit{boundary} $\partial M = M\backslash \operatorname{int} (M)$ of $M$.  Observe that if $p \in \mathcal{V} M$  then for any chart $(U_\alpha, \varphi_\alpha)$  we have that $x_{1}=0$ for $x=\varphi_\alpha^{-1}(p)$ .
This easily follows from the property that an immersion of an open subset of  $\mathbb{R} ^{n}$  into  $\mathbb{R} ^{m}$  takes smooth curves into smooth curves.

Denote by  $C^\infty (M)$ the  $\mathbb{R}$-algebra of smooth functions  $f:M\rightarrow \mathbb{R}$. As for a classical manifold, for  $p\in M$  we define its \textit{tangent space} by
\begin{align*}
T_p M := \left\{v  : C^{\infty}\left(M\right) \rightarrow \mathbb{R} \:|\: v  \text { is a linear map and  }   v(fg)=v(f)g(p)+f(p) v(g) \ \forall f, g \in C^{\infty}\left(M\right)\right\}.
\end{align*}
For a chart $(U_\alpha, \varphi _{\alpha})$ of $M$ and $p\in \varphi _\alpha ( U_\alpha)$, define
$$
\partial_{i}(p): f \in C ^\infty(M) \mapsto \frac{\partial f\circ \varphi _{\alpha} }{\partial x_{i}}\left(\varphi _{\alpha} ^{-1}(p)\right) ,\quad  i=1, \ldots , n. 
$$
Then the set $\{\partial _1 (p), \ldots , \partial _n (p)\}$  is a coordinate frame at  $T_p M$.

\begin{figure}[H]
\centering
\includegraphics[scale=0.19]{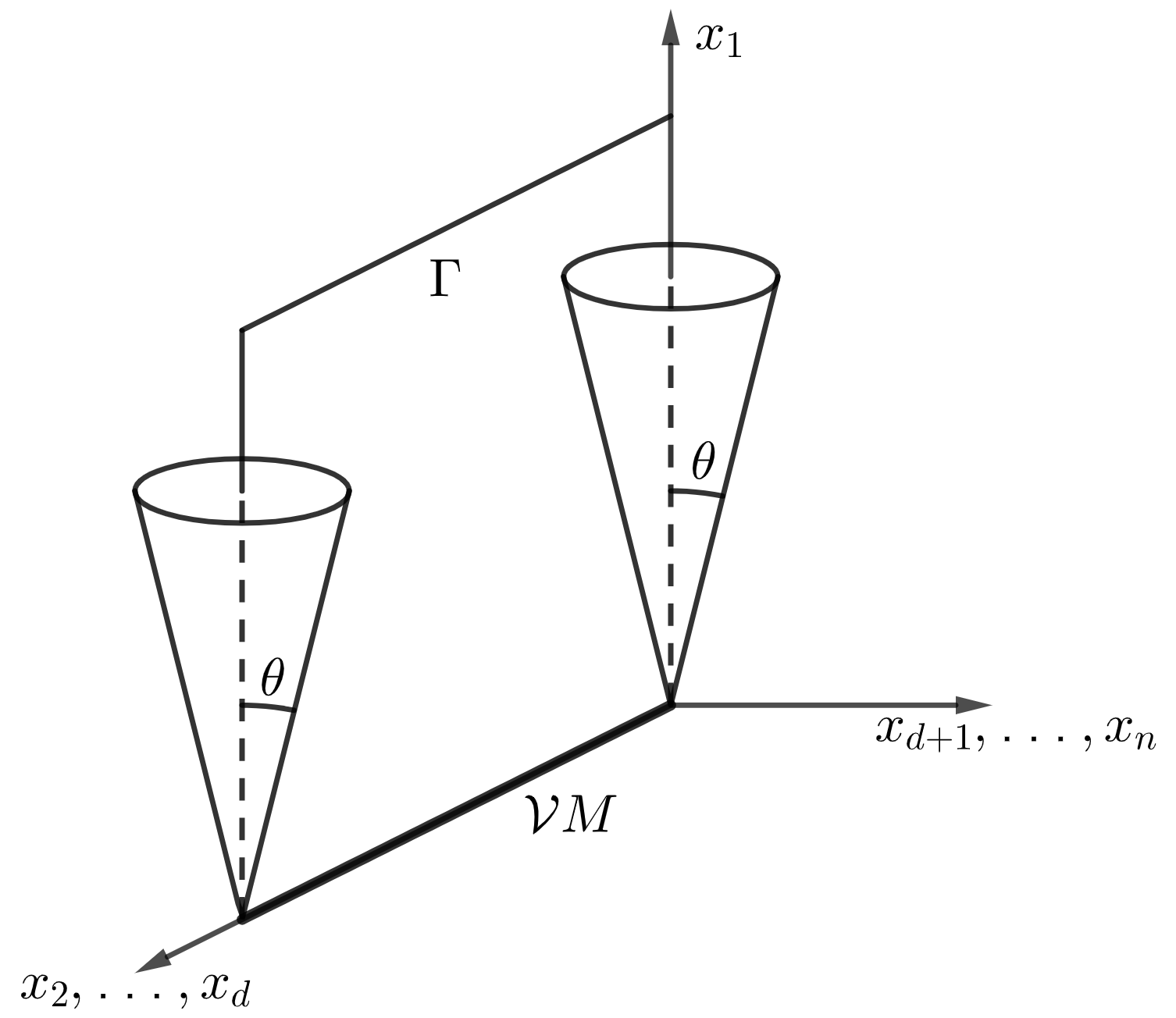}
\caption{A local chart of a s.c.s.s. }
\label{218}
\end{figure}

\begin{definition} \label{107}
Let $M^n$ be a compact smooth solid conical manifold and let $\left\{\Gamma _m ^{d_m} \right\}$ be a finite disjoint family of $d_m$-dimensional  compact submanifolds of $M$ with boundary $\partial \Gamma _m = \Gamma _{m} \cap \partial M $ such that 
\begin{enumerate}[$(i)$]
\item $\cup _{m} \partial \Gamma _{m} = \mathcal{V}M$.
\item If $p\in \partial\Gamma _{m}$ and $\left(U_{\alpha} , \varphi _{\alpha} \right)$ is a chart of $M$ with $p\in \varphi _\alpha (  U_{\alpha} )$, then 
$$\left. \varphi _{\alpha} \right|_{U_{\alpha} \cap \{x_{d_{m}+1}=\cdots = x_n =0\}}$$ 
is a chart of $\Gamma_m$; see Figure \ref{218}.
\end{enumerate}
We define the \textnormal{singular set} of  $M$ as $\Gamma = \cup _{m} \Gamma _{m}$ and call the pair $(M, \Gamma)$ a \textnormal{solid conical singular space (s.c.s.s.)}.  
\end{definition}

A \textit{Riemannian metric} on a smooth solid conical manifold  $M$  is a correspondence which associates to each point  $p$  of  $M$  an inner product  $g( \  , \ )_{p}$ (that is, a symmetric, bilinear, positive definite form) on the tangent space   $T_{p} M$,  which varies smoothly in the following sense:
if  $\left(U_{\alpha},\varphi _{\alpha} \right)$  is a  chart of $M$ at $p$ then  $x\mapsto g\left(\partial  _i , \partial _j  \right)_{\varphi _\alpha (x)}$   is a smooth function on  $U_{\alpha}$.

\begin{definition} \label{107a}
A \textnormal{Riemannian s.c.s.s.} is a s.c.s.s. $(M , \Gamma)$ endowed with a Riemannian metric $g$, satisfying the following  conformal condition along  $\mathcal{V}M$:  for all $m$  and $p\in \partial \Gamma _{m}$ there is $(U_\alpha, \varphi_\alpha)$, where $U_\alpha$ is open in  $\mathcal C_{d_m,h_m}$, such that for some positive function $\varrho$,
\begin{equation} \label{237} 
\varphi ^{\ast} _{\alpha} g(x) = \varrho(x)\delta _{\mathbb{R}^n} \  \text { for all } \  x\in U_{\alpha} \cap \{x_1 = x_{d_m +1}= \cdots x_n =0\},
\end{equation}
where $\delta _{\mathbb{R}^n}$ is the Euclidean metric.
\end{definition}
This definition simply ensures that, $h_m = \cot \theta_m$, where $\theta_m$ is the angle between $\Gamma _m $ and  $\partial M$ along $\mathcal{V} M$, calculated with respect to the metric $g$.

The following technical lemma will be used later:

\begin{lemma} \label{106} 
Let $(M,\Gamma)$ be a s.c.s.s.. Then we have:
\begin{enumerate}[$(i)$]
\item There exists  $\varepsilon _0 >0 $ such that $\rho(\cdot) = \operatorname{dist} _g (\cdot , \Gamma )$ is smooth in $\{ p \in M  \:|\: 0 < \rho (p) <\varepsilon _0 \}$,   
$$
\rho (p) = \operatorname{dist} _g (p , \Gamma  _i ) \quad \text { if }  \operatorname{dist} _g (p , \Gamma  _i ) < \varepsilon _0
$$ 
and
$$
\{ p \in M  \:|\:  \operatorname{dist} _g (p , \Gamma  _i )  < \varepsilon _0 \}\cap  \{ p \in M  \:|\:  \operatorname{dist} _g(p , \Gamma  _j ) <\varepsilon _0\} = \varnothing \quad \text { if  } i\neq j. 
$$

\item   In a coordinate neighborhood $\varphi_{\alpha}(U_{\alpha})$, with  $U_{\alpha} \subset \mathcal{C}_{d_m,h_m}$, satisfying \eqref{237}, the following limits hold:
\begin{equation} \label{107b}
\begin{aligned}
g\left(\nabla \rho , \nu \right)_{p} &\rightarrow \frac{h_m}{\sqrt{1+h_m ^2}} \quad \text {  and } \quad  \rho (p) H_{g}(p) &\rightarrow  \frac{(n-d_m -1)h _m}{\sqrt{1+h_m ^2}},
\end{aligned}
\end{equation}
as $p \rightarrow \mathcal{V}M$ along $\partial M \backslash \mathcal{V}M$, where  $\nu$ is the  outward unit normal vector to the boundary  $\partial M\backslash \mathcal{V} M$.

\item For all $\varepsilon >0$  there exists $\delta \in (0,\varepsilon _0) $ such that
$$
\left| (\rho \Delta \rho )(p) - (n-d_i-1) \right|<\varepsilon \quad \text { if }  0< \operatorname{dist} _g (p , \Gamma  _i ) < \delta.
$$ 
\end{enumerate}

\end{lemma}

\noindent \textit{Proof.} \ Set 
$$
\mathcal{C}_{d_m , h_m} ^{\circ} =\left\{\left(x_{1}, \ldots, x_{n}\right) \in \mathbb{R}^{n} \:|\:  x_{1} > h_m \left(x_{d_m +1}^{2}+\ldots+x_{n}^{2}\right)^{\frac{1}{2}}\right\}.
$$

We first construct local extensions of $M$ across $\mathcal VM$.
If $p\in \mathcal{V}M \cap \Gamma _m$ we can assume that there is a chart $\left( U_{\alpha} , \varphi _{\alpha} \right)$ at $p$,  with $U_\alpha = \mathcal{C} _{d_m , h_m} \cap B_\varepsilon$, $B_{\varepsilon} = \{x\in \mathbb{R} ^n \:|\: |x| <\varepsilon \}$, such that there exist smooth functions $ \widehat{g} _{ij} : B_\varepsilon  \rightarrow \mathbb{R}$ satisfying 
$$
 \widehat{g} _{ij} | _{U_\alpha} = g\left( \partial _i , \partial _j \right) \circ \varphi _{\alpha}
$$
and $[\widehat{g}_{ij}]$ is symmetric and positive definite.
Let $\widehat{U}_{\alpha}$ be the set obtained from $\left( B_\varepsilon  \backslash \mathcal{C}_{d_m , h_m} ^{\circ} \right)\cup \varphi _{\alpha}\left( U_{\alpha} \right)$ by identifying the points $x$ and $\varphi _{\alpha}(x)$ whenever $x\in \partial \mathcal{C} _{d_m , h_m} \cap B_\varepsilon $. Set
\begin{equation} \label{172}
\widehat{\varphi} _{\alpha}(x) := \left\{
\begin{aligned}
&\left[ \varphi _{\alpha} (x) \right] & & \text { if }   x\in U_{\alpha} , \\
&\left[ x \right] & & \text { if }  x\in B_\varepsilon  \backslash \mathcal{C} _{d_m , h_m} ^{\circ}.
\end{aligned}
\right.
\end{equation}
Observe that $\widehat{U} _{\alpha}$ is a manifold without boundary with a Riemannian metric defined by
$$
 \widehat{g} \left(\widehat{\partial} _i , \widehat{\partial} _j \right):=  \widehat{g} _{ij}\circ\widehat{\varphi} ^{-1} _{\alpha}. 
$$ 
where $\left\{\widehat{\partial} _{i}\right\}$ is the coordinate frame associated to $\widehat{\varphi} _{\alpha}$. Similarly, if $\mathbb{R}^{d_m} := \mathbb{R}^n \cap \{x_{d_m +1}=\cdots=x_n =0\}$, we can extend 
$\varphi _{\alpha} (U_{\alpha} \cap \mathbb{R}^{d_m})\subset \Gamma$ to a smooth manifold without boundary denoted by $\widehat{U_{\alpha} \cap \mathbb{R}^{d_m}}$. 

Choosing $\varepsilon $  smaller if necessary, we can prove $(ii)$ and also that $\widehat{\rho} |_{\varphi _{\alpha (U_{\alpha})}} = \rho | _{\varphi _{\alpha (U_{\alpha})}}$. The item $(i)$ easily follows from this.

The proof of $(iii)$ can be found in \cite[pp. 257]{LN}, and concludes the proof of the lemma.

\begin{flushright}
$\blacksquare$
\end{flushright}

We end this section returning to our model, Example \ref{model:ex}, which is given by the compact set $\mathcal{C} _{d,h} \cup \{\infty\}$ endowed with the metric 
$$
g_0 (x):= \rho (x) ^{-2} \delta _{\mathbb{R}^n},
$$ 
where  $\rho (x) = \operatorname{dist} _{\delta _{\mathbb{R}^n} }(x , \Gamma )=\sqrt{x_{d+1} ^2 + \cdots + x_{n} ^2 }$.

Observe that the inversion $x\mapsto x / |x|^2$ provides a chart on a neighborhood of $\infty$ and it is also an isometry of $g_0$.
We finally determine the mean curvature $H_{g_0}$ of $\partial \mathcal{C} _{d,h} \backslash \mathcal V \mathcal{C} _{d,h}$:

\begin{lemma}
We have  
\begin{equation} \label{125}
H_{g_0}=- \frac{dh}{ \sqrt{1 + h^2}}.
\end{equation}
\end{lemma}

\noindent \textit{Proof.} \ Let  $(\mathcal{O},\Theta) $  be a chart of  $\mathbb{S} ^{n-d-1}$  and let  $\{\partial _{\theta _i} \}$  be its associated coordinate frame.  Set $U=\{(r,y,z,s)\in \mathbb{R}\times  \mathbb{R}^{d-1} \times  \mathbb{R}^{n-d-1}\times  \mathbb{R}  \:|\: z\in \mathcal{O}, \  r>0 \text { and } -r/h< s \}$  and define 
$$
\overline{\varphi}: (r , y , z , s) \in U \mapsto \left( hr  - s  ,  y  ,   (h s  +  r )\Theta (z) \right)\in \mathbb{R} ^n.
$$
We see that  $\overline{\varphi}$ is a chart of $\mathbb{R}^n\backslash \mathbb{R}^d$  and  $\varphi (r , y , z) := \overline{\varphi} (r , y , z , 0)$  is a chart of  $\partial \mathcal{C}  _{d,h} \backslash \mathcal{V} \mathcal{C} _{d,h}$. Let  $\{\partial _i\}_{i=1} ^n$  be the associated coordinate frame to  $\overline{\varphi}$.  We can see   $\nu _{g_0} (r ,  y , z) = \left(|\partial _{n}|_{g_0} ^{-1} \partial _{n} \right)(r , y ,z , 0)$ is the outward unit normal vector to the boundary  $\partial \mathcal{C}  _{d,h} \backslash \mathcal{V} \mathcal{C}  _{d,h}$. We have
\begin{align*}
g_0(\partial _1 , \partial _1) (r , y ,z, s) &=  \frac{1 + h^2}{(hs+r)^2}, \quad g_0(\partial _i , \partial _i)(r , y , z, s) =  \frac{1}{(hs+r)^2},\\
g_0(\partial _j , \partial _k) (r , y , z, s) &= \delta _{\mathbb{R}^n} \left( \partial _{\theta _{j-d}} , \partial _{\theta _{k-d}} \right)(z), \quad  g_0(\partial _{n},\partial _{n})(r , y , z , 0) = \frac{1 + h^2}{r^2},
\end{align*}
where  $i = 2 , ... , d$ and $j,k= d+1 , ... , n-1$.  Observe also that $\{\partial _i\} _{i=1}^{d}$ is orthogonal.

Finally, the mean curvature is given by
\begin{equation*}\label{152}
H_{g_0} = \frac{1}{2}\sum _{i, j=1}^{n-1} g_0 ^{ij}  \nu _{g_0} \left( g_{0,ij} \right)= \frac{1}{2} \sum _{i=1}^{d} \frac{r}{\sqrt{1+h^2}} \frac{\partial g_{0,ii}}{\partial s}=  - \frac{dh}{ \sqrt{1 + h^2}} .
\end{equation*}
This proves \eqref{125}.

As a final remark, notice that
$$
H_{\delta_{\mathbb{R}^n}}(r , y , z) = \frac{(n-d-1)h}{\sqrt{1+h^2}} \frac{1}{r}, 
$$
so that the mean curvature, when calculated in terms of $\delta_{\mathbb R^n}$, blows up as we approach the singular set.

\begin{flushright}
$\blacksquare$
\end{flushright}



\section{The main theorem} \label{210}

In this section, we state and prove our main result:
\begin{theorem}\label{209} 
 Let $\left(M,\Gamma,g\right)$ be a Riemannian $s.c.s.s.$ with $n\geq 3$ and $\Gamma = \cup _m \Gamma _m$.  Let $c_0\in L^{\infty} (M\backslash \Gamma)$ and $c_1 \in L^{\infty} ( \partial  M \backslash \mathcal{V} M )$  be smooth functions bounded below by a positive constant. There exists a metric $\tilde{g}$ on $M\backslash \Gamma$,   conformal to $g$ and complete as a metric space, with scalar curvature $R_{\tilde{g}}=-c_0$ in $M\backslash  \Gamma$ and mean curvature $H_{\tilde{g}}=-c_1$ on $\partial M \backslash \mathcal{V}M$ if and only if  $\operatorname{dim} \Gamma _m >\frac{n-2}{2}$ for all $m$.
\end{theorem}

To that end, we will prove the existence of a positive smooth solution  $u$  of the following problem:
\begin{equation}\label{22} 
\left\{ 
\begin{aligned}
-\Delta _g u + \frac{(n-2)}{4(n-1)} R_g u  + c_0 u^{\frac{n+2}{n-2}}  &= 0  & & \text { in }   M\backslash \Gamma  , \\
\frac{\partial u}{\partial \nu _g}  + \frac{n-2}{2(n-1)} H_g  u + c_1 u^{\frac{n}{n-2}} &= 0 & &  \text { on }  \partial M \backslash\mathcal{V}M , \\
\liminf _{p \rightarrow \Gamma} u(p)\operatorname{dist} _g (p, \Gamma)^{\frac{n-2}{2}} &>0.
\end{aligned}
\right.
\end{equation}
Observe that the inequality above ensures that $\widetilde g=u^{\frac{4}{n-2}}g$ is complete as a metric space on $ M\backslash \Gamma$.

Our proof goes along the same lines as \cite{aviles}. We also provide the necessary modifications of some results in \cite{1aviles}. 
For any open set $U\subset M\backslash \Gamma$ such that $U \cap \partial M \neq \varnothing$, we set
$$
\mathcal{D} _{U}= \partial U \cap \operatorname{int} (M) \quad \text { and } \quad  \mathcal{N}_{U} = U\cap \partial M 
$$
and define
\begin{equation}\label{75}
 \lambda_g \left(U\right)= \inf_{\zeta \in C_c ^\infty \left(U \right)}\frac{\int_{U} \left(|\nabla \zeta|^2 _g+\frac{n-2}{4(n-1)}R_g \zeta ^2\right)\dv _g + \frac{n-2}{2(n-1)}\int _{\mathcal{N}_{U}} H_g \zeta ^2 \s _g}{\int_{U} \zeta^2 \dv_g},
\end{equation}
where $\dv _g$ and $\s _g$ are the volume and the area elements respectively. We divide the proof of the existence part in two cases, according to the sign of $\lambda_g \left(M\backslash \Gamma \right)$. The non-negative and negative cases are handled in Subsections \ref{158} and \ref{159}, respectively. In Subsection \ref{onlyif} we prove the converse statement, i.e., that the existence of a positive solution to \eqref{22} implies that  $\operatorname{dim} \Gamma _m >\frac{n-2}{2}$ for all $m$.


\subsection{The case $\lambda _g \left(M \backslash \Gamma \right) \geq 0 $}  \label{158}
Unless otherwise stated, $\Omega$ will denote an open subset of $M\backslash \Gamma$ satisfying the following: 
\begin{assumption} \label{214}The closure $\overline{\Omega}\subset M\backslash \Gamma$ is a   compact $n$-submanifold  with corners locally modeled on  $\mathbb{R}^{n-2} \times [0,\infty)^2$ (see the Appendix \ref{105} for the precise definition) and  $\overline{\mathcal{D}_{\Omega}}$ and $\overline{\mathcal{N}_{\Omega}}$  are  $(n-1)$-submanifolds with smooth boundaries  $\partial \overline{\mathcal{D}_{\Omega}} := \overline{\mathcal{D} _{\Omega}} \backslash \mathcal{D} _{\Omega}$ and $\partial \overline{\mathcal{N}_{\Omega}} := \overline{\mathcal{N} _{\Omega}} \backslash \mathcal{N} _{\Omega}$, respectively (observe that $\partial \overline{\mathcal{D}_{\Omega}} = \partial \overline{\mathcal{N}_{\Omega}} \subset \partial M \backslash \mathcal{V} M$).
\end{assumption}

Recall that $\rho (\cdot )= \operatorname{dist} _g (\cdot , \Gamma)$ and set
$$
M_j=\left\{ p\in M \:|\: \rho (p) > j^{-1} \right\}, \quad j=1, 2, \ldots.
$$
Since $\Gamma$ is compact and $0\leq g(\nabla \rho , \nu ) <1 $ near $\mathcal{V}M$ (see the Lemma \ref{106} $(ii)$), for large  $j$, $M_j$   is an open set in $M$ satisfying the Assumption \ref{214}.

\begin{figure}[H]
\centering
\includegraphics[scale=0.16]{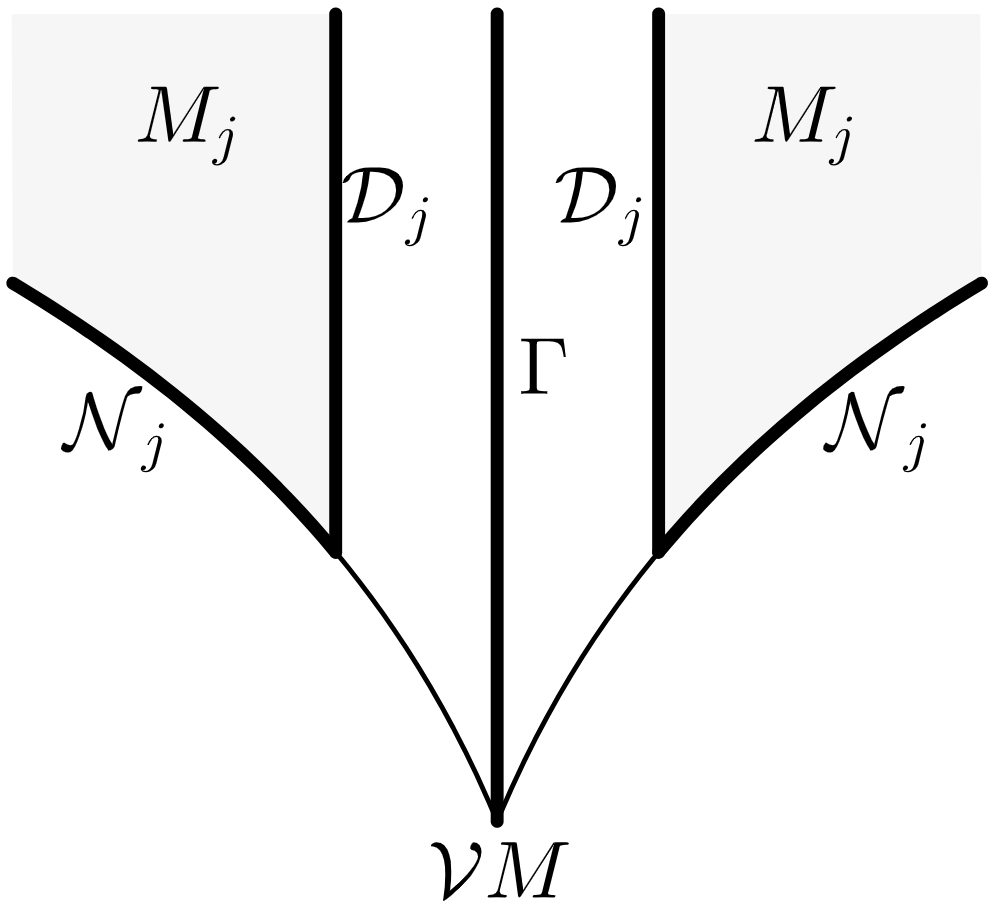}
\caption{ }
\label{252}
\end{figure}

For simplicity we write,
\begin{gather*}
\mathcal{D}_j = \mathcal{D} _{M_j},  \quad  \mathcal{N}_j =\mathcal{N}_{M_j}, \quad \mathcal{D} = \mathcal{D} _{\Omega} \quad \text { and } \quad \mathcal{N} =\mathcal{N}_{\Omega};
\end{gather*}
see Figure \ref{252}. The following maximum principle is proved as in \cite[Chapter 2, Section 5]{prot} (see also \cite{rosales}):

\begin{proposition}  \label{maxlem2}
Let $U$ be an open bounded set of $\mathbb{R}^n$ satisfying $\partial U = \overline{\mathcal{D}} \cup \overline{\mathcal{N}}$, $\mathcal{D}\cap \mathcal{N}  = \varnothing$  and $\partial \mathcal{D}= \partial \mathcal{N}$,  where  $\mathcal{N}$  is a hypersurface (without boundary) in  $\mathbb{R}^{n}$, and $\mathcal{D}\neq \varnothing$.  Let $\{a^{ij}\}$ be symmetric and uniformly elliptic coefficients with $ a^{i j} \in C^{0, 1-\frac{n}{q}}(U)$,  $i, j \in\{1, \ldots, n\}$, $q>n$.
Suppose that $d\in L^\infty (U)$, $\hat{c}_0 \in L^2 (U)$,  $\hat{c}_0 \geq 0$,  $\hat{c}_1,  \hat{c}_2 : \mathcal{N} \rightarrow \mathbb{R}$,  $\hat{c}_1\geq 0$. 
Consider the operators
$$
Lu:= -\left(	a^{ij} u_{x_i} \right) _{x_j}   - du
\quad \text { and } \quad
Bu:=\frac{\partial u}{\partial \nu} + \hat{c}_2 u  ,
$$
and let $f_0 , f_1 :\mathbb{R}\rightarrow \mathbb{R}$ be non decreasing functions.
Assume that
$$u, v , w \in C^1 \left(U\cup \mathcal{N}\right)\cap C \left(\overline{U}\right),$$
$f_0(u), f_0(v) \in L^2(U)$,  and 
\begin{equation*}
	\left\{
	\begin{aligned}
	Lu  + \hat{c}_0f_0 (u) 	&\leq Lv+ \hat{c}_0f_0(v) &  & \text { in } U ,\\
	u&\leq v & & \text { on }  \mathcal{D} , \\	
	 Bu + \hat{c}_1 f_1(u) &\leq Bv +\hat{c}_1 f_1(v) &  & \text { on }  \mathcal{N},
 	\end{aligned}
	\right.
	\end{equation*}
	\begin{equation*}
	\left\{
	\begin{aligned}
	Lw & \geq 0 &   & \text { in } U ,\\
	Bw &\geq 0 & & \text { on }    \mathcal{N} ,\\
	w&>0	& & \text { in }   \overline{U},
	\end{aligned}
	\right.
	\end{equation*}	
are satisfied in the weak sense in $U$, and in the classical sense on $\mathcal{N}$. Then
$$
u\leq v \quad  \text { in }  U.
$$		
\end{proposition}

Set
$$
\mathcal{L} _g u :=-\Delta _g u + \frac{(n-2)}{4(n-1)} R_g u
\quad \text { and } \quad
\mathcal{B} _g u := \frac{\partial u}{\partial \nu _g}  + \frac{n-2}{2(n-1)} H_g u .
$$

\begin{lemma} \label{A}
If  $u, v\in C^2 \left(\Omega\right)\cap  C \left(\overline{\Omega}\right) $  are non-negative functions such that 
	\begin{equation*}
	\left\{
	\begin{aligned}
	\mathcal{L} _g u + c_0 u^{\frac{n+2}{n-2}} &\leq \mathcal{L} _g v + c_0 v^{\frac{n+2}{n-2}}   & & \text { in }  \Omega ,\\
	u&\leq v    & & \text { on }  \mathcal{D} ,\\
	\mathcal{B} _g u +  c_1 u^{\frac{n}{n-2}} & \leq \mathcal{B} _g v + c_1v^{\frac{n}{n-2}} & & \text { on }   \mathcal{N},
	\end{aligned}
	\right.
	\end{equation*}
then
$$
u\leq v.
$$
\end{lemma}

\noindent \textit{Proof.}  
Choose $j$ large enough such that $\overline{\Omega} \subset M_j$. The condition  $\lambda _g (M \backslash \Gamma ) \geq 0 $ implies $\lambda _g (M_j) \geq 0$. Using standard variational arguments (see for example the proof of Proposition \ref{propo:exist} in the appendix), we can find $\phi_j \in C^{\infty} (M_j ) \cap C\left(\overline{M_j}\right)$ a minimum of the variational problem
\begin{equation*}
 \lambda_{g} \left(M_{j}\right):= \inf_{\zeta \in C_c ^\infty \left(M_{j} \right)}\frac{\int_{M_{j}} \left(|\nabla \zeta|^2 _{g}+\frac{n-2}{4(n-1)}R_{g} \zeta ^2\right) \dv _g + \frac{n-2}{2(n-1)}\int _{\mathcal{N}_{j}} H_{g} \zeta ^2 \s _g }{\int_{M_{j}} \zeta^2 \dv _g },
\end{equation*}
so that 
\begin{equation} \label{137}
\left\{
\begin{aligned}
\mathcal{L} _g \phi _j &= \lambda _g (M_j) \phi _j & & \text {  in }  M _j   , \\
\phi _j &= 0 & & \text { on }  \mathcal{D}_j  ,\\
\mathcal{B} _g \phi _j &= 0  &  & \text { on }    \mathcal{N} _j ,\\
\phi _j & > 0  & & \text { in }    M_{j} .
\end{aligned}
\right.
\end{equation}
In particular, the conditions of Proposition \ref{maxlem2} are fulfilled with $w=\phi_j$ and the result follows.

\begin{flushright}
$\blacksquare$
\end{flushright}

\begin{lemma} \label{B}
There exists  $u \in C^\infty\left(\Omega\right)$  such that
\begin{equation} \label{103} 
\left\{ 
\begin{aligned}
\mathcal{L} _g u  + c_0 u^{\frac{n+2}{n-2}}   &= 0 & &   \text { in }  \Omega ,\\
u(x)&\rightarrow \infty                             &  & \text { as }  x  \to \overline{\mathcal{D}} , \\
\mathcal{B} _g u +c_1 u^{\frac{n}{n-2}}  &= 0       &  & \text { on }   \mathcal{N} ,\\
u&>0 & & \text { in }  \Omega.
\end{aligned}
\right.
\end{equation}
\end{lemma}
\noindent \textit{Proof.}   
As in the last proof, we choose $j$ large such that  $\overline{\Omega} \subset  M_j$ and $\phi _{j+1} $ satisfying \eqref{137}.  Then  $\widetilde{g} = \phi _{j+1} ^{\frac{4}{n-2}} g $ satisfies
$$
\left\{
\begin{aligned}
R_{\widetilde{g}}  = \frac{4(n-1)}{(n-2)} \lambda _g (M_{j+1}) \phi _{j+1} ^{-\frac{4}{n-2}} &\geq 0 & & \text {  in   }  M _{j} ,\\
H_{\widetilde{g} } = \phi _{j+1} ^{- \frac{n}{n-2}} \left( \frac{2(n-1)}{n-2}   \frac{\partial \phi _{j+1}}{\partial \nu _g} + \phi _{j+1}   H_{g}\right)  &= 0 & & \text { on }  \mathcal{N} _j.
\end{aligned}
\right.
$$
It follows from Proposition \ref{propo:exist} that, for every positive integer  $m$,  there is $w_m \in C^{2,\alpha} \left(\Omega \right) \cap C^{0,\alpha}\left(\overline{\Omega}\right)$  such that 
$$
\left\{
\begin{aligned}
 \mathcal{L} _{\widetilde{g}}w_m+ c_0 w_m ^{\frac{n+2}{n-2}} &=0 & &  \text { in }   \Omega ,\\
w_m &= m & & \text { on }  \overline{ \mathcal{D}} ,\\
\mathcal{B} _{\widetilde{g}} w_m + c_1 w_m ^{\frac{n}{n-2}} & =0 &  &\text { on }   \mathcal{N} , \\
w_m &>0 & & \text { in }   \overline{\Omega}.
\end{aligned}
\right.
$$
By the Lemma \ref{A}, the  $w_m$ are monotonic increasing and by Theorem  \ref{14} the $w_m$ are uniformly bounded on any open subset $V\subset \Omega$ with $\overline{V} \subset \Omega$.
Using elliptic estimates as in the last paragraph of the proof of Proposition \ref{propo:exist}, we can assume that $w_m$ converges $C^{2}$ in compact subsets of $\Omega$ to a solution $w$ of 
\begin{equation} \label{145}
\left\{
\begin{aligned}
\mathcal{L}  _{\widetilde{g}}  w+  c_0 w ^{\frac{n+2}{n-2}}  &=0 &  & \text { in } \Omega ,\\
w(x) &\rightarrow \infty  & & \text { as }   x\rightarrow \overline{\mathcal{D}} ,\\
\mathcal{B}  _{\widetilde{g}}  w + c_1 w ^{\frac{n}{n-2}}&=0 & & \text { on }   \mathcal{N} ,\\
w&>0 & & \text { in }   \Omega.
\end{aligned}
\right.
\end{equation}
Hence,  $u=\phi _{j+1} w$  solves \eqref{103}, and the regularity of $u$ follows from standard elliptic arguments.

\begin{flushright}
$ \blacksquare $
\end{flushright}
 
It follows from Lemma \ref{B} that, for each $j$, we can find $u_j \in C^\infty (M _j)$ satisfying
\begin{equation*} 
\left\{ 
\begin{aligned}
\mathcal{L} _g u_j + c_0 u_j ^{\frac{n+2}{n-2}}  &=0 & & \text { in }  M _j ,\\
u _j (x) &\rightarrow \infty                          & & \text { as }  x  \to \overline{\mathcal{D} _j} , \\
\mathcal{B} _g u_j + c_1 u_j ^{\frac{n}{n-2}}  &= 0  & & \text { on }   \mathcal{N} _j , \\
u_j &>0 & & \text { in }  M_j.
\end{aligned}
\right.
\end{equation*}
By Lemma \ref{A}, the $u_j$ are monotonically decreasing, i.e., $u_j \geq u_{j+1}$ in $M_j$. Similarly as we did in the proof of Lemma \ref{B}, for $u(x) := \inf_{m\geq j} u_m (x)$, $x\in M_j$, we have that $u\in C^\infty \left( M\backslash \Gamma\right)$ and
\begin{equation*}
\left\{
\begin{aligned}
 \mathcal{L} _gu +  c_0u^{\frac{n+2}{n-2}} &=0   & &  \text { in }  M\backslash \Gamma ,\\
\mathcal{B} _g u + c_1 u ^{\frac{n}{n-2}}&=0     & & \text { on }  \partial M\backslash \mathcal{V} M ,\\
u& \geq 0                                        & & \text { in }  M\backslash \Gamma.\\
\end{aligned}
\right.
\end{equation*}

It remains prove that $u>0$ and $\liminf _{p\rightarrow \Gamma} u(p)\rho ^{\frac{n-2}{2}} (p)>0$. To that end, we set
$$
\psi := C_{\ast} \rho ^{\frac{2-n}{2}},
$$
where the constant $C_{\ast}>0$ will be determined below. By the Lemma \ref{106},  there exists $\varepsilon >0$ small enough such that 
\begin{equation}
\rho \Delta _g \rho-\frac{n}{2}+\frac{1}{n-1} \|R _g\|_{L^\infty} \rho^{2} <-C_1  \label{165}
\end{equation}
and
\begin{equation}
-  g\left(\nabla \rho , \nu \right) + \frac{1}{n-1}\rho H_g < -C_1, \label{166}
\end{equation}
hold on $\left\{ p\in M \:|\: 0< \rho (p) < \varepsilon \right\}$, where $C_1$ is a positive constant independent of $C_{\ast}$. Indeed, the inequality \eqref{165} is a consequence of 
$\min _k \{ \operatorname{dim} ( \Gamma _k ) \}>\frac{n-2}{2} $ and Lemma \ref{106}$(iii)$, that implies $(\rho (p)\Delta \rho (p)-\frac{n}{2}) \rightarrow (\frac{n-2}{2} - d_j)$ as $p\rightarrow \Gamma _j$. The inequality \eqref{166} is a consequence  of Lemma \ref{106}$(ii)$.

On  $M_{m} \backslash \overline{ M_{m_0}}=\left\{ p\in M \:|\: m^{-1} < \rho (p) < m_0 ^{-1} \right\}$, $\varepsilon ^{-1}<m_0 < m$, we have
\begin{align*}
&\mathcal{L} _g \psi + c_0 \psi ^{\frac{n+2}{n-2}}  - \frac{n-2}{4(n-1)} C_{\ast}  m_0 ^{-\frac{2-n}{2}} R_g   \\
 & \leq C_{\ast} \frac{(n-2)}{2} \rho^{-\frac{n+2}{2}}\left(\rho \Delta _g \rho-\frac{n}{2}\right) + \|c_0\|_{L^{\infty}} C_{\ast} ^{\frac{n+2}{n-2}} \rho ^{-\frac{n+2}{2}} + \frac{n-2}{2(n-1)}   C_{\ast} \rho^{\frac{2-n}{2}}\|R_g \|_{L^\infty}
\\ 
&\leq  C_{\ast} \rho ^{-\frac{n+2}{2}}\left[\frac{n-2}{2}\left(\rho \Delta _g \rho-\frac{n}{2}+\frac{1}{n-1} \|R_g \|_{L^\infty} \rho^{2}\right) + \|c_0\|_{L^{\infty}} C_{\ast} ^{\frac{4}{n-2}}  \right]\\
&\leq C_{\ast} \rho ^{-\frac{n+2}{2}}\left[\frac{n-2}{2}\left(-C_1+\frac{1}{n-1} \|R_g \|_{L^\infty} \rho^{2}\right) + \|c_0\|_{L^{\infty}} C_{\ast} ^{\frac{4}{n-2}}  \right], 
\end{align*}
by \eqref{165}. Moreover,
\begin{align*}
&\mathcal{B} _g \psi + c_1 \psi ^{\frac{n}{n-2}}  -  \frac{n-2}{2(n-1)}C_{\ast} m_0 ^{-\frac{2-n}{2}}H_g \\
   &\leq   \frac{n-2}{2} C_{\ast} \left(-  g\left(\nabla \rho , \nu \right) + \frac{1}{n-1}\rho H_g  +\|c_1\| _{L^\infty}\frac{2}{n-2} C_{\ast} ^{\frac{2}{n-2}}  \right)  \rho ^{-\frac{n}{2}} 
   + \frac{(n-2)}{2(n-1)}C_{\ast} m_0 ^{-\frac{2-n}{2}} \rho^{-1} o(1) \\
 &\leq  C_{\ast} \rho ^{-\frac{n}{2}} \frac{n-2}{2}\left( -C_1  +     \|c_1\| _{L^\infty}  C_{\ast} ^{\frac{2}{n-2}}  + (m_0 \rho)^{\frac{n-2}{2}} o(1)\right),
\end{align*}
above we have used   Lemma \ref{106}$(ii)$  and  \eqref{166}. Then, for $m_0$ large enough and  very small $C_{\ast}$ we have that 
$$
\left\{
\begin{aligned}
\mathcal{L} _g \psi   + c_0 \psi ^{\frac{n+2}{n-2}} - \frac{n-2}{4(n-1)} C_{\ast} m_0 ^{-\frac{2-n}{2}} R_g  &\leq 0 & &  \text { in }   M_m \backslash \overline{ M_{m_0} },\\
\mathcal{B} _g \psi + c_1 \psi ^{\frac{n}{n-2}}  - \frac{n-2}{2(n-1)} C_{\ast}  m_0 ^{-\frac{2-n}{2}} H_g & \leq 0  & & \text { on }  \mathcal{N}_m \backslash \overline{ \mathcal{N}_{m_0} }. 
\end{aligned}
\right.
$$
Applying the Lemma \ref{A},
$$
u_{m} \geq \psi-C_{\ast} m_0 ^{-\frac{2-n}{2}} \quad \text { on }  M_m \backslash \overline{ M_{m_0} }.
$$
In particular,
$$
u\geq \psi-C_{\ast}  m_0 ^{-\frac{2-n}{2}} \quad \text { on }  \left\{p\in M\:|\: 0< \rho (p)\leq m_0 ^{-1} \right\}.
$$
Therefore,  
$$
\liminf _{p\rightarrow \Gamma} u(p)\rho ^{\frac{n-2}{2}} (p)>0,
$$
and, by maximum principle arguments, $u>0$ on $M\backslash \Gamma$.


\subsection{The case  $\lambda _g \left(M\backslash \Gamma  \right) < 0 $}  \label{159}
 Let $\varepsilon _0 >0 $ be given by the Lemma \ref{106}$(i)$, and let $m_0 > \varepsilon _0 ^{-1}$ be such that
 \begin{equation*}
\rho ^2 R_g  + 2(n-1) \left(\rho \Delta_g \rho - \frac{n}{2}\right) <-C_2 
\end{equation*}
and
\begin{equation*}
\rho H_g  -(n-1) g\left(\nabla \rho , \nu \right) <-C_2, 
\end{equation*}
on $\{p\in M \:|\: 0< \rho (p) <m_0 ^{-1}\}$, where $C_2$ is a positive constant.

Let $\zeta \in C^{\infty} _c \left(M \backslash \overline{M_{m_0 +1 }}\right)$  be such that
\begin{equation*}
\left\{
\begin{array}{cl}
0\leq \zeta \leq 1 &    \text { in } M ,\\
\zeta  = 1   &   \text { on } M \backslash \overline{M_{m_0 +2 }}.
\end{array}
\right.
\end{equation*}
Set
\begin{equation*} 
\widetilde{\rho} (p)=\left\{\begin{array}{rl}
\left(\rho (p) - 1 \right)\zeta (p)  +1 , &   \text { if }  p\in M \backslash  \left( \overline{M_{m_0 +1 }} \cup \Gamma\right),\\
1 , &   \text { if }   p \in M_{m_0 + 1 },
\end{array}\right.
\end{equation*}
and $\widetilde{g} :=  \widetilde{\rho} ^{-2} g$. On  $M \backslash  \left( \overline{M_{m_0 +2 }} \cup \Gamma\right)$   we have
\begin{equation*}
R_{\widetilde{g}}  = \rho ^2 R_g  + 2(n-1) \left(\rho \Delta_g \rho - \frac{n}{2}\right)
\end{equation*}
and
\begin{equation*}
\begin{aligned}
H_{\widetilde{g}} =\rho H_g  -(n-1) g\left(\nabla \rho , \nu \right).
\end{aligned}
\end{equation*}
By Lemma \ref{106}$(ii)$, observe that $H_{\widetilde{g}} \in L^\infty$. For small $\varepsilon >0$ we have
$$
\left\{
\begin{aligned}
R_{\widetilde{g}}\varepsilon + c_0 \varepsilon ^{\frac{n+2}{n-2}} & \leq 0 & & \text { in } M \backslash  \left( \overline{M_{m_0 +2 }} \cup \Gamma\right) ,\\
H_{\widetilde{g}}\varepsilon + c_1 \varepsilon ^{\frac{n}{n-2}}&\leq 0 & & \text { on } \partial M \backslash  \left( \overline{\mathcal{N}_{m_0 +2 }} \cup \mathcal{V} M\right).
\end{aligned}
\right.
$$

\subsubsection*{Step 1} \label{205}
It is easy to see that the property $\lambda _g (M \backslash \Gamma) <0$ is preserved under conformal changes of the metric. Choose $j_1 > m_0 + 2 $ such that $\lambda _g (M_{j_1}) <0$. As in the positive case, we choose $\phi _{j_1 }\in C^{2,\alpha} (M_{j_1} ) \cap C^{0,\alpha}\left(\overline{M_{j_1}}\right)$ to be a minimum of the variational problem
\begin{equation*}
 \hat{\lambda}_{\widetilde{g}} \left(M_{j_1}\right):= \inf_{\zeta \in C_c ^\infty \left(M_{j_1} \right)}\frac{\int_{M_{j_1}} \left(|\nabla \zeta|^2 _{\widetilde{g}}+\frac{n-2}{4(n-1)}R_{\widetilde{g}} \zeta ^2\right) \dv _{\widetilde{g}} + \frac{n-2}{2(n-1)}\int _{\mathcal{N}_{j_1}} H_{\widetilde{g}} \zeta ^2 \s _{\widetilde{g}} }{\int_{M_{j_1}} \zeta^2 \dv  _{\widetilde{g}}  + \int _{\mathcal{N}_{j_1}} \zeta ^2 \s  _{\widetilde{g}}}.
\end{equation*}
Then  $\hat{\lambda}_{\widetilde{g}} \left(M_{j_1}\right) <0$ and 
$$
\left\{
\begin{aligned}
\mathcal{L} _{\widetilde{g}}\phi _{j_1}&=\hat{\lambda}_{\widetilde{g}} \left(M_{j_1}\right)  \phi_{j_1}  & & \text { in }  M_{j_1}  ,\\
\phi_{j_1} &=0 & & \text { on }    \mathcal{D}_{j_1} , \\
\mathcal{B} _{\widetilde{g}} \phi _{j_1} & = \hat{\lambda}_{\widetilde{g}} \left(M_{j_1}\right) \phi _{j_1} &  & \text { on }   \mathcal{N}_{j_1}, \\
\phi_{j_1}&>0 & & \text { in }   M_{j_1}.
\end{aligned}
\right.
$$ 
Choose  $\delta >0$ such that 
$$
\left\{
\begin{aligned}
\|c_0\|_{L^{\infty}}\left(\delta \phi _{j_1} \right)^{\frac{4}{n-2}} &\leq -\hat{\lambda}_{\widetilde{g}} \left(M_{j_1}\right) & & \text { in }  M_{j_1}  , \\
 \|c_1\|_{L^{\infty}} \left(\delta \phi _{j_1} \right)^{\frac{2}{n-2}} &\leq -\hat{\lambda}_{\widetilde{g}} \left(M_{j_1}\right) & &  \text { on } \mathcal{N} _{j_1} . 
\end{aligned}
\right.
$$
If  $\psi _{j_1} := \delta \phi _{j_1}$,
$$
\left\{
\begin{aligned}
\mathcal{L} _{\widetilde{g}} \psi _{j_1} + c_0 \psi _{j_1} ^{\frac{n+2}{n-2}}&\leq 0 &  &\text { in  }  M _{j_1} ,\\
\psi_{j_1}&=0 & & \text { on }   \mathcal{D}_{j_1} ,\\
\mathcal{B}  _{\widetilde{g}} \psi_{j_1} + c_1 \psi _{j_1} ^{\frac{n}{n-2}} & \leq 0 &  & \text { on }   \mathcal{N}_{j_1} , \\
\psi_{j_1}&>0 & & \text { in }     M_{j_1} .
\end{aligned}
\right.
$$

Now we use monotone iteration schemes as in \cite[Section 2.3]{sattin}. For a large positive number $S$ there exist $u_m ^{+}, u_m ^{-} \in  C^{2,\alpha} ( M _{j_1} ) \cap C^{0,\alpha}\left(\overline{M_{j_1}}\right)$ such that
\begin{equation}\label{89}
\psi _{j_1}=u_0 ^- \leq u_1 ^- \leq \cdots \leq u_m ^- \leq u_m ^+ \leq \cdots \leq u_1 ^+ \leq u_0 ^+ = S
\end{equation}
and
\begin{equation}\label{86}
 \left\{
 \begin{aligned} 
-\Delta _{\widetilde{g}}  u_{m+1} ^{\pm}  + S^{\frac{n+2}{n-2}}  u_{m+1} ^{\pm}   &= F_m   & & \text { in }  M_{j_1} ,\\
u_{m+1} ^{\pm}  &=\varepsilon & & \text { on }    \mathcal{D}_{j_1}  ,\\
\frac{\partial u_{m+1} ^{\pm} }{\partial \nu _{\widetilde{g}}} + S^{\frac{n}{n-2}} u_{m+1} ^{\pm} &= G_m  & & \text { on }    \mathcal{N}_{j_1},
\end{aligned}
\right. 
\end{equation}
where $u_{0} ^-  = \psi _{j_1}$, $u_{0} ^+ =S$,
\begin{gather*} 
F_{m} := \left(S ^{\frac{n+2}{n-2}}- R_{\widetilde{g}} \right) u_{m}^{\pm} - c_0\left|u_{m}^{\pm}\right|^{\frac{n+2}{n-2}} \quad \text { and } \quad
 G_{ m } := \left(S ^{\frac{n}{n-2}}- H _{\widetilde{g}}\right)u_{ m}^{\pm} - c_1\left|u_{m}^{\pm}\right|^{\frac{n}{n-2}}. 
\end{gather*}
Indeed, we can assume that 
\begin{gather*}
t\in \left[0,S\right] \mapsto \left(S ^{\frac{n+2}{n-2}} - R _{\widetilde{g}}(p) \right)t - c_0 (p) t^{\frac{n+2}{n-2}} \quad \text { and } \quad
t\in \left[0,S\right] \mapsto \left(S ^{\frac{n}{n-2}} - H _{\widetilde{g}}(q) \right)t - c_1 (q) t^{\frac{n}{n-2}}
\end{gather*}
are non-decreasing functions.  So, the inequalities \eqref{89} are consequences of the following of the maximum principle, which is proved using \cite[Theorem 8.16]{gilb}) and \cite[Hopf boundary point lemma]{rosales}:

\begin{proposition} \label{147}
Assume the hypotheses of Proposition \ref{maxlem2}, and $d\in L^\infty (U)$ with $d\leq 0$. Let $\hat{c}_2:\mathcal{N} \rightarrow [0,\infty) $ be a function and  $u\in H^1 (U)\cap C(\overline{U})$. Suppose 
$$
\lim _{t\rightarrow 0^-} \frac{u(x + t \nu (x)) - u(x)}{t}
$$ 
exists for all $x\in \mathcal{N}$ and 
\begin{equation*}
	\left\{
	\begin{aligned}
-\left(	a^{ij} u_{x_i} \right) _{x_j}    - du   &\leq   0  &  & \text { in } U ,\\
u&\leq 0  &   & \text { on }  \mathcal{D} , \\	
	 \lim _{t\rightarrow 0^-} \frac{u\left(x + t \nu(x)\right) - u(x)}{t} + \hat{c}_2 (x) u(x) &\leq  0 & & \text { on }  \mathcal{N},
 	\end{aligned}
	\right.
	\end{equation*}
is satisfied in the weak sense in $U$, and in the classical sense on $\mathcal{N}$, where  $\nu$  is the outward unit normal vector to $\mathcal{N}$. Then  
$$
u \leq 0 \quad \text { in }  U.
$$
\end{proposition}

By \eqref{89}, both sequences $\{u _m ^+\}$ and $\{u _m ^-\}$  converge pointwise, so we set
$$
u^- := \lim _{m\rightarrow \infty} u_m ^- \quad \text { and } \quad u^+ := \lim _{m\rightarrow \infty} u_m ^+.
$$
Using standard elliptic estimates, we see that the sequences $\left\{u^{\pm} _m\right\}$ are uniformly bounded in $C^{2,\alpha} (M _{j_1} ) \cap C^{0,\alpha}\left(\overline{M_{j_1}}\right)$. It follows that $u^{\pm}\in C^{2,\alpha'} ( M _{j_1} ) \cap C^{0,\alpha'}\left(\overline{M_{j_1}}\right)$ for $\alpha'<\alpha$ and
$$
\left\{
\begin{aligned} 
\mathcal{L} _{\widetilde{g}} u^{\pm}+c_{0} \left(u^{\pm}\right)^{\frac{n+2}{n-2}} &=0 &  & \text { in }  M_{j_1} ,  \\
 u^{\pm} &= \varepsilon &  & \text { on }  \mathcal{D}_{j_1} ,\\
\mathcal{B} _{\widetilde{g}} u^{\pm}  + c_1 \left( u ^{\pm } \right)^{ \frac{n}{n-2}} &= 0 &  & \text { on }  \mathcal{N}_{j_1} ,\\
u^{\pm}  &>0 & &  \text { in } \overline{ M_{j_1} }.
\end{aligned}
\right.
$$

\subsubsection*{Step 2}

For $k> j_1$, we define
\begin{equation*}
u_{k0} ^{-} (p) = \left\{
\begin{array}{rl} 
u^{-} (p) , & \text { if }   p\in M_{j_1} ,\\
\varepsilon ,&  \text { if }   p\in M_{k}\backslash M_{j_1} . 
\end{array}\right.
\end{equation*}
Proceeding as in Step 1 (replacing $\psi _{j_1}$ by $u^- _{k0}$),  there exist  $u^{\pm} _{k} \in C^{2,\alpha} \left( M _{k}\right)\cap C^{0,\alpha}\left(\overline{M _{k}}\right)$  such that
$$
\left\{
\begin{aligned} 
\mathcal{L}_{\widetilde{g}}  u^{\pm} _k +c_{0} \left(u^{\pm} _k\right)^{\frac{n+2}{n-2}} &=0 & & \text { in }  M _k  , \\
u^{\pm} _k &\geq \varepsilon & & \text { in } M_k \backslash \overline{M_{j_1}} , \\
\mathcal{B} _{\widetilde{g}} u^{\pm} _k  + c_1 \left( u ^{\pm } _k \right)^{ \frac{n}{n-2}} &= 0 &  & \text { on }  \mathcal{N}_{k} , \\
u^{\pm}  _k &>0 &  & \text { in } \overline{M_{k}}.
\end{aligned}
\right.
$$
 
Using Theorem  \ref{14} and standard elliptic arguments, we see that there is $u\in C^{\infty} (M\backslash \Gamma)$ such that
\begin{equation*}
\left\{
\begin{aligned} 
 \mathcal{L} _{\widetilde{g}} u +c_{0} u^{\frac{n+2}{n-2}}  &=0 &  & \text { in }  M\backslash \Gamma , \\
u&\geq \varepsilon &  & \text { in }   M\backslash \left( \overline{M_{j_1}} \cup \Gamma\right) , \\
 \mathcal{B} _{\widetilde{g}} u+c_{1} u^{\frac{n}{n-2}} &=0 & &  \text { on }  \partial M \backslash \mathcal{V} M , \\
u&>0 & & \text { in }  M\backslash \Gamma.
\end{aligned}
\right.
\end{equation*}

Recalling that we are working with the metric $\widetilde{g} =  \widetilde{\rho}^ {-2} g$, we have that $w= u\widetilde{\rho}^{\frac{2-n}{2}} $ satisfies
\begin{equation*}
\left\{ 
\begin{aligned}
-\Delta _g w + \frac{(n-2)}{4(n-1)} R_g w  + c_0 w^{\frac{n+2}{n-2}}  &= 0  & & \text { in }   M\backslash \Gamma ,  \\
\frac{\partial w}{\partial \nu _g}  + \frac{n-2}{2(n-1)} H_g w + c_1 w^{\frac{n}{n-2}} &= 0 & &  \text { on }   \partial M\backslash\mathcal{V}M  , \\
\liminf_{p \rightarrow \Gamma} w(p)\operatorname{dist}_g (p, \Gamma)^{\frac{n-2}{2}} &\geq \varepsilon.
\end{aligned}
\right.
\end{equation*}
This completes the proof of the case $\lambda_g(M\backslash \Gamma)<0$.


\subsection{The $``$only if$"$ part}\label{onlyif}
Following  the lines of \cite{aviles} we suppose that $d_r \leq (n-2)/2$, for some $r$, and $\widetilde{g}= u^{\frac{4}{n-2}}g$ has scalar curvature $R_{\widetilde{g}} =-c_0$ on $M\backslash \Gamma$ and mean curvature $H_{\widetilde{g}} =-c_1$ on $\partial M \backslash \mathcal{V}M$. 
Then $u$ satisfies 
\begin{equation} \label{243}
\left\{
\begin{aligned}
\Delta _g u &= \frac{n-2}{4(n-1)}\left( R_g +  \frac{c_0}{2} u^{\frac{4}{n-2}}  \right) u  +  \frac{n-2}{4(n-1)}  \frac{c_0}{2} u^{\frac{n+2}{n-2}} & & \text { in }  M\backslash \Gamma ,\\
\frac{\partial u}{\partial \nu _g} &= - \frac{n-2}{2(n-1)} H_g u  - \frac{n-2}{2(n-1)} c_1 u^{\frac{n}{n-2}} & & \text { on } \partial M \backslash \mathcal{V}M
\end{aligned}
\right.
\end{equation}
and
\begin{equation}\label{242}
\begin{split}
& -\int_{M} g(\nabla u, \nabla \zeta)\dv _g =   \int _M  \frac{n-2}{4(n-1)}\left( R_g +  \frac{c_0}{2} u^{\frac{4}{n-2}}  \right) u \zeta \dv _g + \int _M  \frac{n-2}{4(n-1)}  \frac{c_0}{2} u^{\frac{n+2}{n-2}}  \zeta \dv _g \\
&+ \int_{\partial M } \frac{n-2}{2(n-1)} \left(  H_g  +  c_1 u^{\frac{2}{n-2}} \right) u \zeta  \s _g,
\end{split}
\end{equation}
for all $\zeta \in C^1 _c ( M \backslash \Gamma )$.

Assume by contradiction that $\widetilde{g}$ is complete as a metric space. In particular, $u(p)\rightarrow \infty$ as $\rho (p)= \operatorname{dist} _g (p, \Gamma_r)\rightarrow 0$. By Lemma \ref{106}$(ii)$, we can assume that 
\begin{equation}\label{H:pos}
H_g \geq 0 \quad \text { near } \Gamma_r.
\end{equation}

From \eqref{243}, for $\delta >0$ small enough, $u$ satisfies
\begin{equation} \label{244}
\left\{
\begin{aligned}
\Delta _g u &\geq    C_0 u^{\frac{n+2}{n-2}},\\
\frac{\partial u}{\partial \nu _g} &\leq -  C_0 u ^{\frac{n}{n-2}},
\end{aligned}
\right.
\end{equation}
in $\left\{ 0<\rho <\delta \right\}$, where $C_0 = C_0 \left(n, c_0, c_1\right)$ is a positive constant.

Let $(U_{\alpha} , \varphi _{\alpha})$ be  a chart  with the property \eqref{237} of Definition \ref{107a}, satisfying $U_\alpha \subset \mathcal{C}_{d_r , h_r}$. 
Let $V\subset \subset U_\alpha \cap \left\{ 0\leq\rho_{\mathbb R^n} <\delta \right\}$ be an open subset. There exists a small constant $k=k(h_r , V ) >0$ that verifies
\begin{equation}\label{247}
\overline{B\left(x_0 , k \rho _{\mathbb{R}^n}(x_0) \right) \cap \mathcal{C}_{d_r , h_r}}\subset U_\alpha \backslash \{x_{d_{r} +1}= \cdots = x_n=0\}
\end{equation}
for all  $x_0 \in V \backslash \{x_{d_{r}+1}= \cdots = x_n=0\}$, where  $\rho _{\mathbb{R}^n}(x_0)= \operatorname{dist} _{\mathbb{R}^n}(x_0, \{x_{d_{r} +1}= \cdots = x_n=0\})$.

Let $x_0\in V\backslash \{x_{d_{r} +1}= \cdots = x_n=0\}$ be fixed. For a suitable constant $C_1 =$ $C_1 (n , c_0 , c_1, g )>0$ the function
$$
w(x)=\frac{C_1 \left(k \rho_{\mathbb{R}^n}(x_0)\right)^{\frac{n-2}{2}}}{\left[ \left(k \rho_{\mathbb{R}^n}(x_0)\right)^2 -|x-x_0|^2\right] ^{\frac{n-2}{2}}}
$$
satisfies
\begin{equation} \label{248}
\left\{
\begin{aligned}
\Delta _g \left(w \circ \varphi _\alpha ^{-1} \right) &\leq    C_0  \left(w \circ \varphi _\alpha ^{-1} \right)^{\frac{n+2}{n-2}} , \\
\frac{\partial \left(w \circ \varphi _\alpha ^{-1} \right)}{\partial \nu _g} &\geq -  C_0 \left(w \circ \varphi _\alpha ^{-1} \right) ^{\frac{n}{n-2}},
\end{aligned}
\right.
\end{equation}
on $\varphi _\alpha \left( B\left(x_0 , k \rho _{\mathbb{R}^n}(x_0) \right) \cap \mathcal{C}_{d_r , h_r}\right)$. 

We will show that $u\leq w \circ \varphi ^{-1} _\alpha$ on  $\varphi _\alpha \left( B\left(x_0 , k \rho _{\mathbb{R}^n}(x_0) \right) \cap \mathcal{C}_{d_r , h_r}\right)$. Suppose by contradiction this inequality does not hold. Define 
$$
B^+ = \left\{p\in \varphi _\alpha \left( B\left(x_0 , k \rho _{\mathbb{R}^n}(x_0) \right) \cap \mathcal{C}_{d_r , h_r} \right) \:|\: u(p)> w \circ \varphi _\alpha ^{-1} (p)\right\} \neq \varnothing
$$ 
and observe that $\overline{B^+}\subset \varphi _\alpha \left( B\left(x_0 , k \rho _{\mathbb{R}^n}(x_0) \right) \cap \mathcal{C}_{d_r , h_r}\right)\subset M\backslash \Gamma$. 

We consider the case $B^+\cap \partial M \neq \varnothing$; the other case follows the same way. By \eqref{244} and \eqref{248},
\begin{equation*}
\left\{
\begin{aligned}
- \Delta _g \left( w \circ \varphi _\alpha ^{-1}  - u\right)&\geq  0 & & \text { in  } B^+, \\
w \circ \varphi _\alpha ^{-1} -u &=0 & & \text { on } \partial B^+ \cap \operatorname{int}(M),\\
\frac{\partial \left(w \circ \varphi _\alpha ^{-1} -u\right)}{\partial \nu _g} &\geq 0 & & \text { on  } B^+ \cap \partial M.
\end{aligned}
\right.
\end{equation*}
It follows from the maximum principle that $u \leq  w \circ \varphi _\alpha ^{-1} $ on $B^+$. Thus, we have a contradiction. So, $u\leq w$  on  $\varphi _\alpha \left( B\left(x_0 , k \rho _{\mathbb{R}^n}(x_0) \right) \cap \mathcal{C}_{d_r , h_r}\right)$. 

In particular,
$$
u\left(\varphi_\alpha (x_0) \right) \leq w(x_0)= C_1\left(k \rho _{\mathbb{R}^n}(x_0)\right)^{-\frac{n-2}{2}}.  
$$ 
We can assume that there is an appropriate constant $C_2>0$ such that $\rho \circ \varphi _\alpha \leq C_2 \rho _{\delta _{\mathbb{R}^n}} $  in $U_\alpha$. Then, since $M$  is compact, there is $\delta _1>0$ small enough so that
\begin{equation}\label{250}
u(p)\leq C_3\rho (p) ^{\frac{2-n}{2}} \quad  \text { if }  0 <\rho (p) < \delta_1,
\end{equation}
where  $C_3$ is a positive constant independent of $u$.

Next we follow the argument in \cite[pp.628-630]{0aviles} to show that $u$ is bounded near $\Gamma_r$. This is where the assumption $d_r\leq \frac{n-2}{2}$ will be used. 
By \eqref{242} and \eqref{H:pos}, there is  $\delta >0$ small enough such that $u$ satisfies
\begin{equation}\label{249}
-\int_{M} g(\nabla u, \nabla \zeta) \dv _g \geq C_4(n,c_0)  \int _M u^{\frac{n+2}{n-2}} \zeta \dv _g ,
\end{equation}
for $\zeta \in C^1 _c (\{ 0 <\rho (p) < \delta\})$ with $\zeta \geq 0$.

Let $\kappa:\mathbb{R}\rightarrow \mathbb{R}\in C^{\infty}$ be  increasing, bounded such that $\kappa(t)\equiv 0$ if $t<0$. Additionally, for each positive integer $m$, let $\xi _m :\mathbb{R}\rightarrow \mathbb{R}$ be a smooth function such that
$$
\xi _m (t) = \left\{\begin{aligned} 
1 & & &\text { if } t> \frac{1}{m},\\
0 & & & \text { if } t<\frac{1}{2m},
\end{aligned}
\right.
$$ 
$0\leq \xi _m \leq 1$ and $|\xi^{\prime} _m|\leq C_5 m $, for a suitable positive constant $C_5$. Let $\beta > \max _{ \{\rho = \rho _0\} }u$  where $\rho _0 \in ( 0 , \delta _1)$. By \eqref{249}, if $m> 3 / \rho _0$,

\begin{align}
&C_4 \int _{M} u^{\frac{n+2}{n-2}}\left(\kappa(u-\beta) \xi _m\circ \rho\right) \dv _g \leq -\int _{M } g\left( \nabla u , \nabla \left( \kappa(u-\beta) \xi _m\circ \rho \right)\right)\dv _g \nonumber\\
 &=  -\int _{M}  \kappa ^{\prime} (u-\beta) \xi _m \circ \rho |\nabla (u-\beta)^+|^2 \dv_g -\int _{M }  \kappa (u-\beta) \xi _m ^{\prime} \circ \rho g\left(\nabla(u-\beta)^+ , \nabla \rho \right)  \dv _g\nonumber\\
 &\leq  C_6  m\|\nabla(u-\beta)^+\|_{L^2 \left(\left\{\frac{1}{2m} \leq \rho \leq \frac{1}{m}\right\}\right)}  \left| \left\{\rho \leq \frac{1}{m} \right\}  \right|^{\frac{1}{2}}. \label{235}
\end{align}
Next we will estimate the r.h.s. of \eqref{235}. Let $z_m : \mathbb{R}\rightarrow \mathbb{R}$ be a smooth function such that
$$
z_m (t)=\left\{
\begin{aligned}
& 1 & & \text { if }  t \in \left( \frac{1}{2m} , \frac{1}{m} \right),\\
& 0 & & \text { if } t \in \left(-\infty , \frac{1}{4m}\right) \cup \left( \frac{5}{4m} , \infty \right),
\end{aligned}
\right.
$$
$0\leq z_m \leq 1$ and $|z_m ^\prime| \leq C_7 m$, where $C_7$ is a positive constant. Since $(u-\beta)^{+} \left( z_m \circ \rho  \right)^2 \in H^1 (M)$ is non-negative and has support in $M\backslash \Gamma$, by \eqref{249}, 
$$
\int _M \left( z_m \circ \rho  \right)^2 |\nabla (u-\beta)^{+}|^2 \dv _g+ \int _M 2(u-\beta)^{+} \left(z_m \circ \rho \right) g\left(\nabla (u-\beta)^{+}, z_m ^\prime \nabla  \rho   \right) \dv _g\leq 0.
$$ 
Thus,
\begin{equation}\label{236}
\|\nabla(u-\beta)^+\|_{L^2 \left(\left\{\frac{1}{2m} \leq \rho \leq \frac{1}{m}\right\}\right)}  \leq C_8 m \|(u-\beta)^+\|_{L^2 \left(\left\{\frac{1}{4m} \leq \rho \leq \frac{5}{4m}\right\}\right)}. 
\end{equation}
By \eqref{250}, \eqref{235} and \eqref{236},
$$
\int _{M} u^{\frac{n+2}{n-2}}\kappa(u-\beta) \xi _m\circ \rho \dv _g \leq C_9 m^2  m^{\frac{n-2}{2}}\left| \left\{\rho \leq \frac{5}{4m} \right\}  \right|  \leq C_{10} m^{2+\frac{n-2}{2} + d_r -n}<C_{11}.
$$
where $C_{11}$ is a positive constant independent of $u$ and $m$. So $\int _M u^{\frac{n+2}{n-2}} <\infty$. By \eqref{235} and   \eqref{250},
\begin{align*}
&\int _{M} u^{\frac{n+2}{n-2}} \kappa(u-\beta)  \xi _m\circ \rho \dv _g\\
& \leq  C_{12}  m\left(\int _{\left\{\frac{1}{2m} \leq \rho \leq \frac{1}{m}\right\}}\left[(u-\beta)^{+}\right] ^{ 2-\frac{n+2}{2(n-2)}}\left[(u-\beta)^+\right]^{\frac{n+2}{2(n-2)}} \dv _g \right)^{  \frac{1}{2}}\left| \left\{\rho \leq \frac{1}{m} \right\}  \right|^{\frac{1}{2}} \\
&\leq C_{13} m\left(\int _{\left\{\frac{1}{2m} \leq \rho \leq \frac{1}{m}\right\}}\left[(u-\beta)^{+}\right] ^{ 4-\frac{n+2}{n-2}} \dv _g \right)^{  \frac{1}{4}} \left| \left\{\rho \leq \frac{1}{m} \right\} \right|^{\frac{1}{2}}\\
&\leq C_{14} m m^{ \frac{n-2}{2}\frac{3n-10}{n-2} \frac{1}{4}} \left| \left\{\rho \leq \frac{1}{m} \right\}  \right|^{\frac{1}{4}+\frac{1}{2}} 
\leq C_{15} m^{\frac{3n-2}{8}} m^{\frac{3}{4} (d_r -n)}=  C_{15}m^{\frac{-3n-2+6d_r}{8}} \leq C_{15} m^{-4}.
\end{align*}

Making $m\to\infty$ we see that $\int _{M} u^{\frac{n+2}{n-2}}\kappa(u-\beta) \dv _g=0$. Therefore, $u$ is bounded near $\Gamma_r$ and we have a contradiction. Hence, $\widetilde{g}$ is not complete as a metric space.


\appendix


\textcolor{white}{\section{Appendix}}
\noindent \textbf{\emph{Appendix}}\\
\subsection{Analysis on cornered manifolds} \label{105} 


In this appendix, we will work with compact cornered $n$-manifolds locally modeled by $\mathbb{R}^{n-2}\times [0,\infty)^2$.  More precisely (see \cite{joyce} for a general definition): 
 
\begin{definition}\label{def:corner}
A smooth \textnormal{cornered manifold} of dimension  $n$ is a paracompact Hausdorff topological space  $\overline{\mathcal{M}}$  and a family of homeomorphisms    
$$
\varphi_{\alpha}: U_{\alpha}  \rightarrow \varphi_{\alpha}\left( U_{\alpha} \right) \subset \overline{\mathcal{M}},
$$  
where  $U_{\alpha}\subset \mathbb{R}^{n-2}\times [0,\infty)^2$ is a relative open subset, satisfying the following three conditions:
\begin{enumerate}[$(i)$]
 \item $\bigcup_{\alpha} \varphi_{\alpha}\left(U_{\alpha}\right)=\overline{\mathcal{M}}$.
 \item For any pair   $\alpha$,  $\beta$,   with  $W:=\varphi_{\alpha}\left(U_{\alpha}\right) \cap \varphi_{\beta}\left(U_{\beta}\right) \neq \varnothing$,  the mappings $\varphi_{\beta}^{-1} \circ \varphi_{\alpha}: \varphi_{\alpha}^{-1}(W) \rightarrow$
$\varphi_{\beta}^{-1}(W)$  is smooth.
 \item The family   $\left\{\left(U_{\alpha}, \varphi_{\alpha}\right)\right\}$  is maximal relative to the conditions $(i)$  and $(ii )$.
 \end{enumerate}
\end{definition}
The pair  $\left(U_{\alpha}, \varphi_{\alpha}\right)$  (or the mapping  $\varphi_{\alpha}$)  with  $p \in \varphi_{\alpha}\left(U_{\alpha}\right)$   is called a \textit{chart}  (or \textit{system of coordinates})  of  $\overline{\mathcal{M}}$  at  $p$.   A family  $\left\{\left(U_{\alpha}, \varphi_{\alpha}\right)\right\}$  satisfying $(i)$  and $(i i)$  is called a \textit{smooth structure} on  $\overline{\mathcal{M}}$. 

\noindent We define the \textit{corner}  by
$$
\mathcal{E}\overline{\mathcal{M}} := \left\{ q\in \overline{\mathcal{M}} \:|\: q =\varphi (y,0,0)   \text { for some  chart  }  (U,\varphi) \text { of } \overline{\mathcal{M}} \text { and } y\in \mathbb{R}^{n-2} \right\}, 
$$
the \textit{interior}
$$
\mathcal{M} =\left\{p \in \overline{\mathcal{M}} \:|\: \exists \left(U_\alpha, \varphi_\alpha\right)   \text { such that } p\in \varphi _{\alpha} \left( U_{\alpha}\right)  \text { and } U_\alpha \text { is an open subset of } \mathbb{R}^n \right\}
$$
and the \textit{boundary} $\partial \overline{\mathcal{M}} = \overline{\mathcal{M}}\backslash \mathcal{M}$ of $\overline{\mathcal{M}}$.

\begin{definition} 
A \textnormal{Riemannian metric} on a smooth cornered manifold $\overline{\mathcal{M}}$  is a correspondence which associates to each point  $p$  of  $\overline{\mathcal{M}}$  an inner product  $g( \  , \ )_{p}$ (that is, a symmetric, bilinear, positive-definite form) on the tangent space   $T_{p} \overline{\mathcal{M}}$, which varies smoothly in the following sense: If  $\left(U_{\alpha} , \varphi _{\alpha} \right)$  is a  chart of $\overline{\mathcal{M}}$ at $p$, and $\left\{\partial _i\right\}$ is the coordinate frame associated to $\varphi _{\alpha}$, then  $g\left(\partial  _i , \partial _j  \right)_{\varphi_\alpha (\cdot)}$ is a smooth function on $U_\alpha$. The manifold $\overline{\mathcal{M}}$ endowed with a Riemannian metric $g$ is called a \textnormal{Riemannian cornered manifold}.
\end{definition}

In what follows, we assume that $(\overline{\mathcal{M}},g)$ is a Riemannian cornered manifold, and
$\partial \overline{\mathcal{M}} = \overline{\mathcal{D} }\cup \overline{\mathcal{N}}$, where $\overline{\mathcal{D}}$ and $\overline{\mathcal{N}}$ are $(n-1)$-submanifolds with  smooth boundary  satisfying 
\begin{equation} \label{224}
\partial  \overline{\mathcal{D}} =\partial  \overline{ \mathcal{N}} =\mathcal{E}\overline{\mathcal{M}} \quad \text { and } \quad \overline{\mathcal{D}} \cap \overline{\mathcal{N}} = \mathcal{E}\overline{\mathcal{M}}.
\end{equation}

For  $1 \leq p < \infty$   we define the Sobolev space  $W^{1, p}(\mathcal{M})$  by
$$
W^{1, p}(\mathcal{M}):=\left\{u \in L^p _{\operatorname{loc}}(\mathcal{M}) \:|\: u\in L^p (\mathcal{M}), \  \exists  \nabla   u    \text { and }     \nabla u \in \overrightarrow{L}^{p}(\mathcal{M})\right\},
$$
equipped with the norm
\begin{equation} \label{92}
\|u\|_{W^{ 1, p} (\mathcal{M})}:=\left(\int _\mathcal{M} |u|^p \dv _g +\int_{\mathcal{M}}|\nabla u|^{p} \dv _g\right)^{\frac{1}{p}}.
\end{equation}

Similarly to the Euclidean case (see \cite{evans, Fan} when $\mathcal{M}$ is an open subset of $\mathbb{R}^n$ with Lipschitz boundary)  we have  the following Sobolev embedding theorems:
\begin{theorem} \label{225}
Assume $1 \leq p<n$, and $u \in W^{1, p}(\mathcal{M})$. Then $u \in L^{p^{*}}(\mathcal{M})$, $p^{*}=\frac{n p}{n-p}$, with the estimate
$$
\|u\|_{L^{p^{*}}(\mathcal{M})} \leq C(p, n, \mathcal{M})\|u\|_{W^{1, p}(\mathcal{M})}.
$$
\end{theorem}

\begin{theorem}\label{226}
(Trace theorem). Assume $1 \leq p<\infty$. Then there is a continuous boundary trace embedding
$$
\mathcal{T}: W^{1, p}(\mathcal{M}) \rightarrow L^{\frac{(n-1) p}{n-p}}(\partial \mathcal{M})
$$
such that
$$
\mathcal{T} u=\left.u\right|_{\partial \mathcal{M}} \quad \text { if } u \in W^{1, p}(\mathcal{M}) \cap C(\overline{\mathcal{M}}).
$$
Moreover, for every $q \in\left[1, \frac{(n-1) p}{n-p}\right)$ the trace embedding $W^{1, p}(\mathcal{M}) \rightarrow L^{q}(\partial \mathcal{M})$ is compact.
\end{theorem}

Set
$$
 H^1 _{\mathcal{D}} (\mathcal{M}):= \left\{u\in H^1 (\mathcal{M})\:|\: u|_{\mathcal{D}} = 0\right\}.
$$
We define the bilinear form
$$
\mathcal{B}[u, v]:= \int _{\mathcal{M}}\left( g\left(\nabla u , \nabla v \right)  + c u v \right) \dv _g + \int_{\mathcal{N}} c_2 u v \s _g  \qquad\left(u, v \in H_{\mathcal{D}}^{1}(\mathcal{M})\right),
$$
where    $c \in  L^\infty (\mathcal{M})$  and   $ c_2 \in  L^\infty (\mathcal{N})$, and we set
$$
\mathcal{B}_{\mu _1, \mu _2} [u,v] := \mathcal{B}[u,v] + \mu _1  \int_{\mathcal{M}} uv  \dv _g + \mu _2 \int _{\mathcal{N}} uv \s _g  \quad (u, v \in H^1 _{ \mathcal{D}} (\mathcal{M})),
$$
for  $\mu _1, \mu _2 \in \mathbb{R}$.

\begin{proposition} \label{26}
There is a number  $\gamma  \geq 0$  such that for each 
$$
\mu _1  \geq \gamma, \quad  \mu _2 \geq 0,
$$
and each function
$$
f_1 \in L^2 (\mathcal{M}) \quad \text { and } \quad f_2 \in L^2 (\mathcal{N}),
$$
there exists a unique weak solution  $u\in H^1 _{ \mathcal{D}} (\mathcal{M})$  of the mixed-boundary-value problem:
$$
\mathcal{B}_{\mu _1, \mu _2 }[u ,v] = \int _{\mathcal{M}} f_1 v  \dv _g +  \int _{\mathcal{N}} f_2 v  \s _g \qquad \left(v\in H^1 _{ \mathcal{D}} (\mathcal{M})\right).
$$
\end{proposition}
\noindent  \textit{Proof.} \ Choose  $\gamma \geq 0$ and $\alpha, \beta >0$ such that
$$
|\mathcal{B}[u,v]| \leq \alpha  \Vert u \Vert _{H^1 (\mathcal{M})} \Vert v \Vert _{H^1 (\mathcal{M})}
$$
and
$$
\beta \Vert u \Vert ^2 _{H^1 (\mathcal{M})} \leq \mathcal{B}[u,u] +  \gamma  \Vert u\Vert ^2 _{L^2(\mathcal{M})} + \mu _2 \|u\|^2 _{L^{2} (\mathcal{N})} ,
$$
hold for all  $\mu _2 \geq 0$  and  $u,  v \in H^1 _{ \mathcal{D}} (\mathcal{M})$.
Then $ \mathcal{B}_{\mu _1, \mu _2} [ \ , \  ]$,  with $\mu _1 \geq \gamma$,  satisfies the hypotheses of the Lax-Milgram Theorem. 

Now fix  $f_1\in L^2 (\mathcal{M})$ and $f_2 \in L^2 (\mathcal{N})$.  Set  $Q(v) :=  \int _{\mathcal{M}} f_1 v \dv _g + \int _{\mathcal{N}} f_2 v \s _g$, which is a bounded linear functional on  $H^1 _{ \mathcal{D}} (\mathcal{M})$. By the Lax-Milgram theorem, there exists a unique  $u\in H^1 _{ \mathcal{D}} (\mathcal{M})$  satisfying 
$$
\mathcal{B}_{\mu _1, \mu _2} [u ,v] =  Q(v),
$$
for all  $v\in  H^1 _{ \mathcal{D}} (\mathcal{M})$.
  
\begin{flushright}
$\blacksquare$
\end{flushright}

\subsection{Elliptic estimates}

In this subsection, we adapt the proof of the next proposition to estimate (in a coordinate neighborhood) our solutions on the boundary. Our  main results are Lemmas \ref{135} and \ref{13}. Denote by $B_{r}$ the ball $\{x\in \mathbb{R}^n\:|\: |x|<r\}$. 

\begin{proposition} \label{55} (See  \cite[Theorem 4.1]{han}  ). \ Suppose $n\geq 2$, $a^{ij}  \in L^\infty (B_1)$  and  $c \in L^q (B_1)$  for some  $q>n / 2$  satisfy
\begin{equation}\label{unif:elliptic}
a^{ij}(x)y_i y_j \geq \lambda |y|^2 \text { for any  } x\in B_1, \ y \in \mathbb{R}^n ,
\end{equation}
and
$$
\Vert a^{ij}\Vert _{L^\infty (B_1)} + \Vert c \Vert _{L^q (B_1)}\leq \Lambda ,
$$
for some positive constants  $\lambda$  and  $\Lambda$. Suppose that  $u\in H^1 (B_1)$  is a subsolution in the sense that the inequality
$$
\int _{B_1} a^{ij} u_{x_i} \zeta _{x_j} + cu\zeta \dx \leq \int _{B_1} f\zeta \dx
$$
holds for any  $\zeta \in H^1 _0 \left(B_1\right)$  and  $\zeta \geq 0$  in  $ B_1$.

If  $f\in L^q(B_1)$,  then  $u^+ \in L^\infty _{\operatorname{loc}} (B_1)$.  Moreover, there holds for any  $r \in (0,1)$  and  $p>0$
$$
\stackbin[B_r]{}{\sup} \   u^+ \leq C \left[ \frac{1}{(1-r)^{\frac{n}{p}}} \Vert u^+\Vert _{L^p (B_1)} + \Vert f \Vert _{L^q(B_1)}\right]
$$
where  $C=C(\lambda , \Lambda , p , q ,n ,B_1) > 0$.
\end{proposition}

For $r>0$, we write
$$
K_r  = \left\{ x  \in \mathbb{R}^n \:|\:  \operatorname{max} _{i\in \{1,\ldots, n\}} |x_i| < r  \right\}.
$$
Set
\begin{gather*}
B=K_1 \cap \left\{ x_{n} > 0 , \ x_{n-1} > 0 \right\},\\
\mathcal{D} = K_1 \cap \left\{x_{n-1}>0, \ x_{n} = 0 \right\} \quad \text { and }  \quad \mathcal{N} = K_1 \cap \left\{x_{n-1} = 0, \  x_n>0 \right\}.
\end{gather*}

Let $n\geq 3$ and $V_0, V_1 \subset \mathbb{R}^n$ be  open bounded sets  with  $V_0 \subset \subset V_1$ and   $\zeta \in C _c ^\infty (V_1)$  satisfying $1\geq \zeta  \geq 0$ and  $ \zeta | _{V_0} \equiv 1$. For the next lemma, we assume that there exist $u\in H^1 (B)$   and  constants $\mathcal{S}_{\ast} \geq 0$,  $C_\ast  >0 $  such that 
\begin{equation} \label{E4}
\left\Vert (u-k) ^+ \zeta \right\Vert _{L^{2^\ast} (B)} \leq C_\ast \left\| \nabla \left[ (u-k) ^+\zeta \right] \right\| _{L^2 (B)},  \text { for all }  k\geq \mathcal{S}_{\ast},
\end{equation}
where $2^{\ast} =\frac{2n}{n-2}$. Write $s_1 = \sup _{\mathbb{R}^n}|\nabla \zeta | $,
$$
\mathtt{V}(k,i)= \{x \in V_i \cap B \:|\: u(x) > k \} \quad  
\text { and } \quad \mathtt{V}_{\mathcal{N}} (k,i)= \{x \in V_i \cap \mathcal{N}   \:|\:  u(x) > k \},  
$$
$i=0,1$. For  $k\in \mathbb{R}$,  set  
$$
u_k = (u-k)^+ .
$$

We also assume that $a^{ij}  \in L^\infty (B)$ is uniformly elliptic respect to $\lambda$ (i.e., \eqref{unif:elliptic} holds),  $c\in L^q (B)$  and  $c_1 \in L^{q_1} (\mathcal{N})$  for some  $q$, $q_1 > n-1$, with
$$
\sum _ {i,j} \Vert a^{ij}\Vert _{L^\infty (B)} + \Vert c \Vert _{L^q (B) } + \Vert c_1 \Vert _{L^{q_1} (\mathcal{N})} \leq \Lambda,
$$
for some positive constant $\Lambda$.

\begin{lemma} \label{lcontb0}
Suppose  $f\in L^q (B)$, $f_1 \in L^{q_1} (\mathcal{N})$  and 
\begin{equation}\label{168}
\int _B \left( a^{ij} u_{x_i} v_{x_j} + cuv \right) \dx + \int  _{ \mathcal{N}}c_1 u v \s \leq \int _B f v \dx + \int_{\mathcal{N}} f_1 v \s  \quad \forall k\geq \mathcal{S}_{\ast},  \ v = u_k \zeta ^2 .    
\end{equation}
Set 
$$
\Psi (k,i) := \sqrt{\int _{\mathtt{V}(k,i)} u_k ^2 \dx + \int _{\mathtt{V}_{\mathcal{N}}(k,i)} u_k ^2 \s }, \quad i=0,1.   
$$
There exist  $\varepsilon = \varepsilon (q , q_1 , n ) >0$ and  $N=N(\lambda , \Lambda , n , B , c , c_1) >0$  such that if  
$$
h > k > N \operatorname{max} \left\{ \Vert u^+ \Vert _{L^2 (B)}  ,   \Vert  u^+ \Vert _{L^2 (\mathcal{N})}  ,  \mathcal{S}_{\ast} \right\}, 
$$ 
then 
$$
\Psi  (h,0) \leq C \left[ s_1  \frac{1}{(h-k)^{\varepsilon}}  + \frac{ \Vert f \Vert _{L^q (B)} + \Vert f_1 \Vert _{L^{q_1} (\mathcal{N})} + h }{(h-k)^{1+\varepsilon}}\right] \Psi ^{1 + \varepsilon } (k,1),
$$
where  $C= C \left( \lambda, \Lambda, n , B , c , c_1 , q , q_1 \right)>0$.
\end{lemma}
\noindent \textit{Proof.}  We will follow the proof in \cite[Theorem 4.1]{han}.
\begin{cl} \label{ccontb1}
We have
\begin{equation*}
\int  _B |\nabla (u_h \zeta )|^2 \dx \leq C_1(\lambda ,\Lambda , n) \left( s_1 ^2 \int _{\mathtt{V}(h,1)} u_h ^2  \dx + \int _B a^{ij} u_{x_i} (u_h \zeta ^2)_{x_j} \dx \right).
\end{equation*}
\end{cl} 
\begin{flushright}
$\square$
\end{flushright}

The inequality \eqref{168} and the Claim \ref{ccontb1} implies  
\begin{equation} \label{169}
\int  _B |\nabla (u_h \zeta )|^2 \dx 
\leq C_1 \left(s_1 ^2 \int _{\mathtt{V} (h,1)} u_h ^2 \dx - \int _B cuv \dx - \int  _{ \mathcal{N}}c_1 u v \s + \int _B f v \dx + \int_{\mathcal{N}} f_1 v \s \right), 
\end{equation}
where $v= u_h \zeta ^2$. Next we will estimate the terms on the right side of \eqref{169}.

\begin{cl}\label{ccontb2}
The following estimates hold:\\

\noindent $(i)$
\begin{align*}
&-\int _{\mathcal{N}} c_1  u (u_h \zeta ^2 ) \s   \\ &\leq   C_2 (n,B) \Vert c_1 \Vert _{L^{q_1} (\mathcal{N})}  |\{u_h \zeta \neq 0\} \cap \mathcal{N} | ^ {\frac{1}{n-1} - \frac{1}{q_1}} \int _{B}  |\nabla (u_h \zeta) |^2  \dx 
+  h^2 \Vert c_1 \Vert _{L^{q_1} (\mathcal{N})} |\{ u_h \zeta \neq 0\} \cap \mathcal{N} |^{1-\frac{1}{q_1}}. 
\end{align*}

\noindent  $(ii)$
\begin{gather*}
\int _{\mathcal{N}} f_1 u_h \zeta ^2 \s 
\leq C_3(n,B) \left( \delta ^{-1} |\{ u_h \zeta \neq 0\}\cap \mathcal{N}| ^ {\frac{n}{n-1} - \frac{2}{q_1}}  \Vert f_1 \Vert _{L^{q_1} (\mathcal{N})} ^2  + \delta  \int _B  |\nabla (u_h \zeta) |^2 \dx \right), 
\end{gather*}
for all  $\delta >0$. \\

\noindent  $(iii)$
\begin{gather*}
-\int _B c u u_h \zeta ^2  \dx 
\leq  C_4 (n,B)\Vert c \Vert _{L^q (B)}  |\{u_h \zeta \neq 0\}| ^ {\frac{2}{n} - \frac{1}{q}} \int _B  |\nabla (u_h \zeta) |^2  \dx +  h^2 \Vert c \Vert _{L^q (B)} |\{u_h \zeta \neq 0\}|^{1-\frac{1}{q}}.
\end{gather*}

\noindent  $(iv)$
\begin{gather*}
\int _B f u_h \zeta ^2 \dx 
\leq  C_5 (n,B) \left( \delta ^{-1} |\{u_h \zeta \neq 0\}| ^ {1 + \frac{2}{n} - \frac{2}{q}} \Vert f \Vert _{L^q (B)} ^2  + \delta  \int _B  |\nabla (u_h \zeta) |^2 \dx\right), 
\end{gather*}
for all $\delta >0$.

\end{cl}
\noindent \textit{Proof.} We will start with the following three inequalities that will be use later.  The Hölder's inequality and the Trace theorem imply
\begin{equation} \label{econtb1}
\Vert u_h \zeta \Vert _{L^2 (\mathcal{N})}  \leq  |\{ u_h \zeta \neq 0 \} \cap \mathcal{N} |^{\frac{1}{2(n-1)}}\Vert  u_h \zeta \Vert _{L^{\frac{2(n-1)}{n-2}} (\mathcal{N})}   
 \leq  C_6 (n,B) |\{ u_h \zeta \neq 0 \}\cap \mathcal{N} |^{\frac{1}{2(n-1)}} \Vert u_h \zeta \Vert _{H^1 (B)}.
\end{equation}
On the other hand, by  \eqref{E4},
\begin{equation} \label{econtb2}
 \int _B (u_h \zeta )^2 \dx \leq  |\{u_h \zeta \neq 0 \}|^{\frac{2^{\ast} -2 }{2^{\ast}}} \|u_h \zeta \|^2 _{L^{2^{\ast}} (B)} 
\leq  C_{\ast} ^2 |\{u_h \zeta \neq 0 \}| ^{\frac{2^\ast - 2}{2^\ast}} \int _B |\nabla (u_h \zeta)| ^2 \dx ,
\end{equation}
for $h\geq \mathcal{S}_\ast$, where $\frac{2^\ast - 2}{2^\ast} = \frac{2}{n}$. Further, from \eqref{econtb1} and \eqref{econtb2},
\begin{equation} \label{253} 
\int _{\mathcal{N}} ( u_h \zeta)^2 \s  \leq C_{7} (n,B) |\{ u_h \zeta \neq 0 \}\cap \mathcal{N} |^{\frac{1}{n-1}}   \int _B |\nabla (u_h \zeta)|^2 \dx
\end{equation}

\noindent $(i)$ \ Since $\frac{n-2}{n-1} + \frac{1}{q_1} <1$,
\begin{align*}
&-\int _{\mathcal{N}} c_1  u (u_h \zeta ^2 ) \s =   -\int _{\{ u_h \zeta \neq 0\} \cap \mathcal{N} } c_1 \left(  u_h ^2 + h  u_h \right) \zeta ^2 \s \\
&\leq 2 \int _{\{ u_h \zeta \neq 0 \} \cap \mathcal{N} }   |c_1|  u_h ^2 \zeta ^2 \s + h^2 \int _{\{ u_h \zeta \neq 0 \} \cap \mathcal{N} }  |c_1|\zeta ^2 \s \\
&\leq  2\left( \int _{\mathcal{N} }  |c_1| ^{q_1} \s \right) ^{\frac{1}{q_1}} \left( \int _{\{ u_h \zeta \neq 0 \} \cap \mathcal{N} }   | u_h \zeta |^ {\frac{2(n-1)}{n-2}} \s  \right) ^{\frac{n-2}{n-1}}  \left|\{ u_h \zeta \neq 0 \} \cap \mathcal{N} \right| ^ {1-\frac{n-2}{n-1} - \frac{1}{q_1}} \\
&+ h^2 \Vert c_1 \Vert _{L^{q_1} (\mathcal{N})} |\{ u_h \zeta \neq 0 \} \cap \mathcal{N} |^{1-\frac{1}{q_1}}\\
&\leq   C_8  (n,B) \Vert c_1 \Vert _{L^{q_1} (\mathcal{N})}  |\{ u_h \zeta \neq 0 \} \cap \mathcal{N} | ^ {\frac{1}{n-1} - \frac{1}{q_1}} \left( \int _B (u_h  \zeta )^2 \dx + \int _{B}  |\nabla (u_h \zeta) |^2 \dx \right)  \\
&+ h^2 \Vert c_1 \Vert _{L^{q_1} (\mathcal{N})} \left|\{ u_h \zeta \neq 0 \} \cap \mathcal{N} \right|^{1-\frac{1}{q_1}}, 
\end{align*}
by \eqref{econtb1}. Using \eqref{econtb2}  we conclude the proof of  $(i)$.\\

\noindent $(ii)$ \ Since  $q_1 > n-1$,    $\frac{1}{q_1} + \frac{n-2}{2(n-1)} < 1$  and  $\zeta \leq 1$,
\begin{equation*}
\int _{\mathcal{N}} f_1 u_h \zeta ^2 \s   \leq   \left( \int _{\mathcal{N}} |f_1| ^{q_1} \s \right) ^{\frac{1}{q_1}} \left( \int _{\mathcal{N}}  | u_h \zeta |^ {\frac{2(n-1)}{(n-2)}} \s  \right) ^{\frac{n-2}{2(n-1)}} \left|\{ u_h \zeta \neq 0\}\cap \mathcal{N}\right| ^ {1- \frac{1}{q_1} - \frac{n-2}{2(n-1)}}.
\end{equation*}
By \eqref{econtb1}  and  \eqref{econtb2},
\begin{align*}
&\int _{\mathcal{N}} |f_1| u_h \zeta ^2  \s  \leq  C_9 \Vert f_1 \Vert _{L^{q_1} (\mathcal{N})}   \left( \int _{B}  |\nabla (u_h \zeta) |^2 \dx  \right) ^{\frac{1}{2}} |\{ u_h \zeta \neq 0\}\cap \mathcal{N}| ^ {\frac{n}{2(n-1)} - \frac{1}{q_1}} \\
&\leq  C_9 \delta ^{-1} |\{ u_h \zeta \neq 0\}\cap \mathcal{N}| ^ {\frac{n}{n-1} - \frac{2}{q_1}}  \Vert f_1 \Vert _{L^{q_1} (\mathcal{N})} ^2  +  C_9 \delta  \int _B  |\nabla (u_h \zeta) |^2 \dx .
\end{align*}
This proves the item $(ii)$.
The proofs of $(iii)$ and $(iv)$ are the same as in \cite{han}.

\begin{flushright}
$\square $
\end{flushright}

Let us observe that  $q > \frac{n}{2}$,  $ q_1 > n-1$,  $\{u_h \zeta  \neq 0\} \subset \mathtt{V}(h,1)$,  $\{u_h \zeta  \neq 0\} \cap \mathcal{N} \subset \mathtt{V} _{\mathcal{N}} (h,1)$,  $|\mathtt{V} (h,1)| \leq h^{-1} \int _{\mathtt{V} (h,1)} u^+  $  and  $|\mathtt{V} _{\mathcal{N}} (h,1)| \leq h^{-1} \int _{\mathtt{V} _{\mathcal{N}}(h,1)} u^+ $. By \eqref{169} and Claim \ref{ccontb2} there  exists a constant   $N = N \left( \lambda , \Lambda , n ,  B , c , c_1\right) >0 $ such that if  $h>N  \max \left\{ \Vert u^+ \Vert _{L^2 (B)}  ,   \Vert  u^+ \Vert _{L^2 (\mathcal{N})}  ,  \mathcal{S}_{\ast} \right\} $  then
\begin{equation}
|\mathtt{V}(h,1)|, \ |\mathtt{V} _{\mathcal{N}} (h,1)| \ <1. \label{econtb5}
\end{equation}
and
\begin{equation} \label{econtb4}
\begin{split} 
&\int  _B |\nabla (u_h \zeta )|^2 \dx \leq   C_{10} \left[ s_1 ^2 \int _{\mathtt{V}(h,1)} u_h ^2    \dx +    \left(  \Vert f \Vert _{L^q (B)} ^2 + h^2 \right) |\{u_h \zeta \neq 0\}|^{1-\frac{1}{q}} \right. \\
 &+ \left(  \Vert f_1 \Vert _{L^{q_1} (\mathcal{N})} ^2 + h^2 \right) |\{ u_h \zeta \neq 0\}\cap \mathcal{N}|^{1-\frac{1}{q_1}}\left] , \textcolor{white}{\int _{\mathtt{V}(1)}} \right.
\end{split}
\end{equation}
where $C_{10}=C_{10}( \lambda ,\Lambda , n , B , c , c_1)>0$.

From   \eqref{econtb1} -  \eqref{253} and  \eqref{econtb4}, we have
\begin{equation}\label{econtb6}
\begin{split} 
& \int  _B (u_h \zeta )^2 \dx  \leq    C_{11}   \left[ s_1 ^2 |\{u_h \zeta \neq 0\}|^{\frac{2}{n}}  \left( \int _{\mathtt{V}(h,1)} u_h ^2 \dx   + \int _{\mathtt{V} _{\mathcal{N}}(h,1)}  u_h ^2 \s \right)  \right. \\
 &  + \left(  \Vert f \Vert _{L^q (B)} + \Vert f_1 \Vert _{L^{q_1} (\mathcal{N})} + h \right)^2 
   \left.    \left( |\{u_h \zeta \neq 0\}|^{1+ \frac{2}{n} -\frac{1}{q}} + |\{u_h \zeta \neq 0\}|^{\frac{2}{n}}  |\{ u_h \zeta \neq 0\}\cap \mathcal{N} |^{1-\frac{1}{q_1}} \right) \right].
\end{split}
\end{equation}
and
\begin{equation}\label{econtb7}
\begin{split}
&\int _{\mathcal{N}} ( u_h \zeta)^2  \s \leq   C_{12}   \left[ s_1 ^2 |\{ u_h \zeta \neq 0\} \cap \mathcal{N} |^{\frac{1}{n-1}} \left( \int _{\mathtt{V} (h,1)} u_h ^2  \dx + \int _{\mathtt{V}_{\mathcal{N}}(h,1)} u_h ^2 \s \right)     \right.  \\
&  + \left(  \Vert f \Vert _{L^q (B)} +  \Vert f_1 \Vert _{L^{q_1} (\mathcal{N})}  + h \right)^2 \\
&   \cdot \left( |\{u_h \zeta \neq 0\} |^{1 -\frac{1}{q}} |\{ u_h \zeta \neq 0\} \cap \mathcal{N} |^{\frac{1}{n-1}} +  \left. |\{ u_h \zeta \neq 0\} \cap \mathcal{N} |^{1 + \frac{1}{n-1} -\frac{1}{q_1}} \right) \right],
 \end{split}
\end{equation}
where $C_i = C_i (\lambda ,\Lambda , n ,B , c , c_1)>0$, $i=11, 12$.

On the other hand. Set $\varepsilon = \frac{1}{n-1} - \frac{1}{\min\{q, q_1\}}$, by Young's inequality,
\begin{equation}\label{206}
\begin{aligned}
 |\{u_h \zeta \neq 0\} |^{\frac{2}{n}}  |\{ u_h \zeta \neq 0\}  \cap \mathcal{N} |^{1 -\frac{1}{q_1}} 
 \leq  \frac{1}{C_{13}}  |\{u_h \zeta \neq 0\} |^{1+\varepsilon} +  \frac{C_{13} -1}{C_{13}} |\{ u_h \zeta \neq 0\} \cap \mathcal{N} | ^{\left(1-\frac{1}{q_1}\right)\frac{C_{13}}{C_{13}-1} }     ,
\end{aligned}
\end{equation}
\begin{equation} \label{254}
\begin{aligned}
 |\{u_h \zeta \neq 0\} |^{1 -\frac{1}{q}} & |\{ u_h \zeta \neq 0\} \cap \mathcal{N} |^{\frac{1}{n-1}}  
\leq \frac{C_{14}-1}{C_{14}} |\{ u_h \zeta \neq 0\}  | ^{\left(1-\frac{1}{q}\right)\frac{C_{14}}{C_{14}-1} }  + \frac{1}{C_{14}}  |\{u_h \zeta \neq 0\} \cap \mathcal{N} |^{1+\varepsilon},
\end{aligned}
\end{equation}
where $C_{13}=\frac{n}{2}(1+\varepsilon)$ and $C_{14}= (n-1)(1+\varepsilon)$. Observe that 
$$
\left(1-\frac{1}{q_1}\right)\frac{C_{13}}{C_{13}-1} \geq 1+\varepsilon \ \text { and  } \ \left(1-\frac{1}{q}\right)\frac{C_{14}}{C_{14}-1}\geq 1+\varepsilon.
$$
 
The inequalities \eqref{econtb5}, \eqref{econtb6} - \eqref{254} imply
\begin{equation}\label{151}
\begin{split}
&\Psi \left(h,0 \right)^2 \leq C_{15}\left[s^2 _1 \left(\left|\left\{u_h \zeta \neq 0\right\}\right| +\left|\left\{u_h \zeta \neq 0\right\}\cap \mathcal{N}\right| \right)^{\varepsilon}\Psi \left(h,1\right) ^2 \right.\\
&  + \left(  \Vert f \Vert _{L^q (B)} +  \Vert f_1 \Vert _{L^{q_1} (\mathcal{N})}  + h \right)^2 \left.\left(\left|\left\{u_h \zeta \neq 0\right\}\right| +\left|\left\{u_h \zeta \neq 0\right\}\cap \mathcal{N}\right| \right)^{1+\varepsilon} \right],
\end{split}
\end{equation}
where $C_{15}= C_{15} \left( \lambda, \Lambda, n , B , c , c_1 , q , q_1 \right)>0$. Consider  \eqref{151} and the following claim, which is proved as in \cite{han}:
\begin{cl}
If $h>k$, then
\begin{gather*}
 |\{u_h \zeta  \neq 0\} | \leq  \frac{1}{(h-k)^2}\int _{\mathtt{V}(k,1)} u_k ^2 \dx, \quad  |\{ u_h \zeta \neq 0\} \cap \mathcal{N} | \leq  \frac{1}{(h-k)^2}\int _{\mathtt{V}_{\mathcal{N}}(k,1)} u_k ^2 \s ,\\
\int _{\mathtt{V}(h,1)} u_h ^2 \dx \leq \int _{\mathtt{V} (k,1)} u_k ^2 \dx  \quad   \text { and } \quad \int _{\mathtt{V} _{\mathcal{N}}(h,1)}  u_h ^2  \s \leq \int _{\mathtt{V} _{\mathcal{N}} (k,1)} u_k ^2 \s .
\end{gather*}
\end{cl}
\noindent Therefore,
$$
\Psi (h,0)^2 \leq C_{16}  \left[ s_1 ^2 \frac{1}{(h-k)^{2\varepsilon}}  + \frac{\left(  \Vert f \Vert _{L^q (B)}  +  \Vert f_1 \Vert _{L^{q_1} (\mathcal{N})} + h \right)^2}{(h-k)^{2(1+\varepsilon)}}\right] \Psi ^{2(1 + \varepsilon )} (k,1),
$$
where $C_{16}= C_{16} \left( \lambda, \Lambda, n , B , c , c_1 , q , q_1 \right)>0$. This proves the Lemma \ref{lcontb0}.
\begin{flushright}
$\blacksquare$
\end{flushright}

Observe that if $u\in H^1 (B)$   and  $ u|_{\mathcal{D}} \in L^\infty (\mathcal{D}) $  then  $ (u - k)^+ |_{\mathcal{D}} =0$  $\forall k\geq \Vert  u \Vert _{L^\infty (\mathcal{D})}$.  By Sobolev embedding inequalities  and Lemma \ref{24} we have that
\begin{equation} \label{E5}
\Vert (u-k)^+ \zeta \Vert _{L^{2^\ast} (B)}\leq C_\ast (n,B)\Vert \nabla \left[ (u-k)^+ \zeta \right] \Vert_{L^2 (B)}   
\end{equation}
for all  $k\geq \Vert u \Vert _{L^\infty (\mathcal{D})}$  and  $\zeta \in C^\infty (\mathbb{R}^n )$. Therefore the condition \eqref{E4} is satisfied.

For a set $A\subset \mathbb{R}^n$  and  $t \in \mathbb{R}$, we write
$$
tA := \{tx \in \mathbb{R} ^n \:|\:  x \in A\} .
$$

\begin{lemma} \label{tcontb1}
Suppose that  $u\in H^1 (B)$,  $u|_{\mathcal{D}} \in L^\infty (\mathcal{D})$,  $\mathcal{S}_{\ast} \geq \Vert u \Vert _{ L^\infty (\mathcal{D}) }$, and
\begin{equation*}\label{E6}
\int _B \left( a^{ij} u_{x_i} v_{x_j} + cuv \right)\dx  + \int  _{\mathcal{N} }c_1 u v  \s \leq \int _B f v  \dx  + \int _{\mathcal{N}} f_1 v \s, \quad v=(u-k)^+ \zeta ^2,
\end{equation*}
for all $k\geq \mathcal{S}_{\ast}$ and $\zeta \in C _c ^\infty \left(K_1\right)$ with $\zeta \geq 0 $. Assume $f\in L^q (B)$  and  $f_1 \in L^{q_1} (\mathcal{N})$,  then if  $p,  p_1 \geq 2$ 
\begin{equation*}
\stackbin[2^{-1}B]{}{\sup} \ u^+  + \stackbin[2^{-1} \mathcal{N}]{}{\sup} \  u^+  \leq 
  C\left( \Vert u^+\Vert _{L^p (B)} +  \Vert  u^+\Vert _{L^{p_1} (\mathcal{N})} +  \mathcal{S}_{\ast} + \Vert f_1 \Vert _{L^{q_1} (\mathcal{N})}  + \Vert f \Vert _{L^q (B)}  \right),  
\end{equation*}
where  $C=C(\lambda , \Lambda , p , p_1 , q , q_1 , n ,  B)>0$.

\end{lemma}
\noindent \textit{Proof.} \ We again follow the lines of  \cite[Theorem 4.1]{han}. As observed above, the condition \eqref{E4} is satisfied. Then, by Lemma \ref{lcontb0} we have  $u^+ \in L^\infty \left(2^{-1}B\right)\cap L^\infty \left(2^{-1}\mathcal{N}\right)$  and
\begin{equation*} 
\stackbin[2^{-1} B]{}{\sup} \ u^+ + \stackbin[2^{-1}\mathcal{N}]{}{\sup} \ u^+  \leq  C \left(  \Vert u^+\Vert _{L^2 (B)} + \Vert  u^+\Vert _{L^2 (\mathcal{N})} + \mathcal{S}_{\ast}  +  \Vert f_1 \Vert _{L^{q_1} (\mathcal{N})}  + \Vert f \Vert _{L^q (B)} \right), \label{E7} 
\end{equation*} 
where $C=C(\lambda , \Lambda , q , q_1 , n ,B)>0$. Using the Hölder's inequality we can conclude the proof.

\begin{flushright}
$\blacksquare$
\end{flushright}

The proof  of the next two lemmas are similar to the one of Lemma  \ref{tcontb1}.

\begin{lemma} \label{135}
Suppose that  $u\in H^1 (B)$,  $u|_{\mathcal{D}} \in L^\infty (\mathcal{D})$,  $\mathcal{S}_{\ast} \geq \Vert u \Vert _{ L^\infty }$, and
\begin{equation*}
\int _B \left( a^{ij} u_{x_i} v_{x_j} + cuv \right) \dx  \leq \int _B f v \dx, \quad v=(u-k)^+ \zeta ^2 ,
\end{equation*}
for all $k\geq \mathcal{S}_{\ast}$ and  $\zeta \in C _c ^\infty \left(K_1\cap \{x_{n-1} >0\}\right)$ and $\zeta \geq 0 $. Assume  $f\in L^q (B)$, then if $p \geq 2$ and $X$ is a compact set with $X\subset B\cup \mathcal{D}$ we have that
\begin{equation*}
\stackbin[X\cap B]{}{\sup}  \ u^+   \leq C\left( \Vert u^+\Vert _{L^p (B)} + \mathcal{S}_{\ast}  + \Vert f \Vert _{L^q (B)}  \right), 
\end{equation*}
where $C=C\left(\lambda , \Lambda , p  , q  , n ,X, B\right)>0$ .
\end{lemma}

We have if  $\zeta \in C_c ^{\infty}\left( K_{1}\cap \{x_n>0\} \right)$,  then  $\zeta |_{\mathcal{D}} = 0$  and
$$
\left\|(u-k)^{+} \zeta \right\|_{L^{2^\ast }(B)} \leq C\left\|\nabla \left[(u-k)^{+} \zeta\right]\right\|_{L^{2}(B)} \quad \forall k \geq 0,
$$
where   $ C=C(n, B)>0$, which implies that the condition  \eqref{E4} is satisfied.

The following lemma will be important in the proof of Theorem \ref{14} below.

\begin{lemma} \label{13}  
Suppose that  $u\in H^1 (B)$  and
\begin{equation*}
\int _B \left( a^{ij} u_{x_i} v_{x_j} + cuv \right) \dx + \int  _{\mathcal{N}}c_1 u v \s \leq \int _B f v \dx + \int _{\mathcal{N}} f_1 v\s , \quad v=(u-k)^+ \zeta ^2, 
\end{equation*}
for all  $k\geq 0$ and  $\zeta \in C _c ^\infty \left(K_1 \cap \{x_n >0\}\right)$  with  $\zeta \geq 0$. Assume  $f\in L^q (B)$   and  $f_1 \in L^{q_1} (\mathcal{N})$, then if  $p,  p_1 \geq 2$ and $X$ is a compact set with $X\subset B\cup \mathcal{N}$ we have that
\begin{equation*}
\stackbin[X\cap B]{}{\sup}  \ u^+ + \stackbin[X \cap \mathcal{N}]{}{\sup} \  u^+  \leq C\left( \Vert u^+\Vert _{L^p (B)} + \Vert  u^+\Vert _{L^{p_1} (\mathcal{N})} + \Vert f_1 \Vert _{L^{q_1} (\mathcal{N})}  + \Vert f \Vert _{L^q (B)}  \right), 
\end{equation*}
where  $C=C\left(\lambda , \Lambda , p , p_1 , q , q_1 , n ,X, B\right)>0$.
\end{lemma}




\subsection{Nonlinear solutions and estimates} \label{19}

Recall that $\overline{\mathcal{M}}=\mathcal{M}\cup\overline{\mathcal D}\cup \overline{\mathcal N}$ denotes a cornered manifold satisfying \eqref{224}. Let   $c$, $c_0 \in L^\infty (\mathcal{M})$, $c_1 \in L^\infty (\mathcal{N})$ be smooth functions  such that $c,  c_0,  c_1 \geq 0$. For $f \in  H^1 (\mathcal{M})$, set
$$
H^1 _{f} (\mathcal{M})= \{ u\in H^{1} (\mathcal{M}) \:|\:   \left. u \right|_{\mathcal{D}} = \left. f \right|_{\mathcal{D}}  \}. 
$$

Our first goal now is to prove the following:
\begin{proposition}\label{propo:exist}
Suppose that $f|_\mathcal{D}\in L^\infty(\mathcal D)$ is non-negative, non-trivial, and H\"older continuous.
Then there exists a solution $u\in C^{2,\alpha}(\mathcal{M}\cup\mathcal N)\cap C(\overline{\mathcal{M}})$ of the mixed boundary problem
\begin{equation} \label{1388}
\left\{
\begin{aligned}
-\Delta u +c u + c_0 u^{\frac{n+2}{n-2}}&=0 &  & \text { in }  \mathcal{M}, \\
u &= f &  & \text { on }  \mathcal{D}, \\ 
\frac{\partial u}{\partial \nu}  + c_1 u^{\frac{n}{n-2}} &=0 & & \text { on }  \mathcal{N},
\end{aligned}
\right.
\end{equation}
satisfying
\begin{equation} \label{142}
u > 0 \quad \text { in }  \mathcal{M} \cup \mathcal{N}. 
\end{equation}
\end{proposition}

We use the following definitions: Let $X$ be a Banach space, a set  $S\subset X$  is said to be \textit{weakly closed} if  $\{u_m\} \subset S$,  $u_m \rightharpoonup u$  implies $u \in S$.  A functional  $\mathcal{I} : S \subset X \rightarrow \mathbb{R}$  is \textit{weakly continuous} at  $u_{0} \in S$  if for every sequence  $\left\{u_{m}\right\} \subset S$  with $u_{m} \rightharpoonup u_{0}$  it follows that  $\mathcal{I}  \left(u_{m}\right) \rightarrow \mathcal{I} \left(u_{0}\right)$.   A functional  $\mathcal{I}  : S\subset X \rightarrow \mathbb{R}$  is \textit{weakly lower semicontinuous} (\textit{w.l.s.c}) at  $u_0 \in S$  if for every sequence  $\{u_m\} \subset S$  for which  $u_m \rightharpoonup u_0$  it follows that  $\mathcal{I}  (u_0) \leq \operatorname{ lim} \operatorname{inf} _{m\rightarrow \infty} \mathcal{I}  (u_m)$.   A functional  $\mathcal{I}  : S\subset X \rightarrow \mathbb{R}$ is \textit{weakly coercive} on  $S$ if  $\mathcal{I}  (u) \rightarrow \infty$  as  $\Vert u \Vert \rightarrow \infty $  on  $S$.  

It is well known that any  closed convex set in a Banach space is weakly closed; see \cite[Theorem 1.39]{yan}.  In particular, $H^1 _{f} (\mathcal{M})$  is weakly closed as it is convex and satisfies
$$
\Vert  u \Vert _{L^{2}(\mathcal{D})}  \leq C\Vert u \Vert _{H^{1}(\mathcal{M})} \quad  \forall u \in H^1 (\mathcal{M}).
$$

\noindent Also, the following Poincaré inequality holds:
\begin{lemma} \label{24}
There exists   $C=C(n,\mathcal{M})>0$  such that
$$
\Vert u \Vert _{L^2 (\mathcal{M})} \leq C \Vert \nabla u  \Vert _{L^2 (\mathcal{M})},
$$      
for all  $u\in H^1 _{\mathcal{D}} (\mathcal{M}):=\left\{v \in H^{1}(\mathcal{M})\:|\: v|_{\mathcal{D}}=0\right\}$.
\end{lemma}
\noindent \textit{Proof.} \ Suppose by contradiction that  there exists a sequence  $\{u_m\} \subset H^1 _{\mathcal{D}} (\mathcal{M})$  such that 
$$
\displaystyle \int _{\mathcal{M}} u_m ^2  \dv _g =1 \quad \text { and } \quad \int _\mathcal{M} |\nabla u_m|^2 \dv _g \rightarrow 0 \quad \text { as }    m\rightarrow \infty.
$$
Hence we may assume  $u_m \rightarrow u_0 \in H^1 (\mathcal{M})$  strongly in  $L^2 (\mathcal{M})$  and weakly in  $H^1 (\mathcal{M})$.  Thus,
$$
\nabla u_0 = 0 \quad \text { and } \quad \Vert u_0 \Vert _{L^2 (\mathcal{M})} =1.
$$
As $H^1 _{\mathcal{D}} (\mathcal{M})$  is convex and closed (the latter being a consequence of Theorem \ref{226}), then
$$
u_0 \in H^1 _{\mathcal{D}} (\mathcal{M}).
$$  
Therefore, $u_0 = 0 $.  This is a contradiction. 
\begin{flushright}
$\blacksquare$
\end{flushright}

Set
$$
F(t):= \frac{n-2}{2n}\left(t^{2}\right)^{\frac{n}{n-2}}  \quad \text { and }  \quad G(t):= \frac{n-2}{2(n-1)}\left(t^{2}\right)^{\frac{n-1}{n-2}}.
$$
Observe that $F$,  $G$ are convex functions and  $F(u)\in L^1 (\mathcal{M})$, $G(u) \in L^1 (\partial \mathcal{M})$  for all $u\in H^1 (\mathcal{M})$. We define the functional  $\mathcal{I}: H^1 _{f} (\mathcal{M}) \rightarrow \mathbb{R}$, by
$$
\mathcal{I} (u) =  \frac{1}{2} \left( \int _{\mathcal{M}} |\nabla u |^2 \dv _g + \int _{\mathcal{M}} c u^{2}  \dv _g \right)+ \int _{\mathcal{M}}c_{0}F(u)  \dv _g + \int _{\mathcal{N}}  c_{1}G ( u) \s _g.
$$
Then $u\in H^1 _{f} (\mathcal{M})$ is a critical point of $\mathcal I$ if and only if it satisfies 
\begin{equation}\label{219}
\int_{\mathcal{M}}\left( g\left(\nabla u , \nabla v \right)+c u v+c_{0}|u|^{\frac{4}{n-2}} u v\right) \dv _g+\int_{\mathcal{N}} c_{1}|u|^{\frac{2}{n-2}} uv \s _g = 0 \quad \forall v \in H_{\mathcal{D}}^{1}(\mathcal{M}),
\end{equation}
i.e., $u$ is a weak solution of \eqref{1388}.
\begin{lemma}\label{propo:I}
The functional $\mathcal{I}$ is $w.l.s.c.$ and weakly coercive.
\end{lemma}

\noindent \textit{Proof.} \ 
By Lemma \ref{24}, 
\begin{equation} \label{220}
C_1 \Vert u \Vert _{H^1(\mathcal{M})} ^2 - C_2 \|f\| _{H^1 (\mathcal{M})} ^2 \leq  \mathcal{I} (u) \quad \forall u\in H^1 _f (\mathcal{M}),
\end{equation}
and further, by Theorems \ref{225} and \ref{226}, 
\begin{align*}
 \mathcal{I} (u) \leq   C_3 \left( \Vert u \Vert _{H^{1} (\mathcal{M})}^{2} + \int _{\mathcal{M}} F(u) \dv _g +  \int _{\mathcal{N}}G( u) \s _g \right) 
\leq C_4 \left( \Vert u\Vert _{H^{1} (\mathcal{M})}^{2} + \Vert u \Vert _{H^{1} (\mathcal{M})}^{2^{\ast}} + \Vert u \Vert _{H^{1} (\mathcal{M})}^{\frac{2(n-1)}{n-2}} \right) ,
\end{align*}
where $C_i =C_i  \left(n,\mathcal{M}\right)$,  $i=1,2$, $C_j = C_j \left( c , c_{0}, c_{1},n,\mathcal{M}\right)$, $j=3 ,4$, are positive constants.  By \eqref{220},  $\mathcal{I}$  is weakly coercive. That $\mathcal{I}$  is  $w.l.s.c.$  follows from the following result:

\begin{proposition}\cite[Theorem 1.41]{yan} \label{28}
Consider the functional  $\mathcal{I}  : C\subset X \rightarrow \mathbb{R}$,  where   $X$  is a real Banach space. Suppose $C$  is closed and convex, $\mathcal{I} $  is convex and continuous. Then   $\mathcal{I} $  is \textit{w.l.s.c.}. 
\end{proposition}

This completes the proof of Lemma \ref{propo:I}.
\begin{flushright}
$\blacksquare$
\end{flushright}

\noindent \textit{Proof of Proposition \ref{propo:exist}.} \ 
Set 
\begin{equation} \label{227}
\ell = \inf _{u\in H^1 _{f} (\mathcal{M})} \mathcal{I} (u).
\end{equation}
It follows from \eqref{220} that $\ell > -\infty$. So, there exists a sequence $\{u_m\} \subset H^1 _{f} (\mathcal{M})$ such that
\begin{equation}\label{228}
\mathcal{I} (u_m)\rightarrow \ell.
\end{equation}
Hence, using \eqref{220} again,  $\{u_m\}$ is uniformly bounded in $H^1(\mathcal{M})$ and so we may assume that this sequence converges weakly in $H^{1}(\mathcal{M})$ to some $u_{0} \in H^{1}(\mathcal{M})$.  Since $ H^1 _{f} (\mathcal{M})$ is weakly closed,  $u_0\in H^1 _{f} (\mathcal{M})$. Then
\begin{equation}\label{229}
\mathcal{I} (u_0) \leq \ell,
\end{equation}
because $\mathcal{I}$ is $w.l.s.c$. So,
\begin{equation}\label{230}
\ell = \mathcal{I} (u_0).
\end{equation}  
Therefore, we  conclude that there exists $u\in H^1 _{f} (\mathcal{M})$ such that \eqref{219} holds.

Now assume that
$$
f |_{\mathcal{D}}\geq 0,
$$
i.e., $f^-:=-\min\{f, 0\}$ satisfies $f^-|_{\mathcal D}=0$.
Observe that $\left|\left( v|_{\mathcal D} \right)\right|= (|v|)|_{\mathcal D}$ for any $v\in W^{1,p} (\mathcal{M}) $.
Since $u^- |_{\mathcal{D}} = -\min \left\{f ,0\right\} = 0$,
$$
\mathcal{I} (u^+  + u^-) = \mathcal{I} (u^+ - u^-)  \quad \text { and } \quad (u^+ + u^-)|_{\mathcal{D}} = f,
$$
we can assume that
$$
u\geq 0.
$$
In particular, by \eqref{219},
\begin{equation} 
\int _{\mathcal{M}} g\left( \nabla u ,\nabla v\right) \dv _g \leq 0 \label{E44} 
\end{equation}
for all  $v\in H^1 _{\mathcal{D}}(\mathcal{M})$  with  $v \geq 0$. Hence, Proposition \ref{55} together with \eqref{E44} gives us $u\in L^\infty _{\operatorname{loc}} (\mathcal{M})$.
Then, the  Lemmas \ref{tcontb1} - \ref{13} imply
\begin{equation}\label{208}
u\in L^{\infty}(\mathcal{M}).
\end{equation}

By \cite[Proposition 2.4]{dron}, $u$ is H\"older continuous on $\overline{\mathcal{M}}$.
Using \cite[Corollary 4.23 and Theorem 4.24]{han} and \cite[Theorem 2]{lieber} we see that $u\in C^{1,\alpha}(\mathcal{M}\cup\mathcal N)$.
Hence, it follows from standard elliptic estimates that $u\in C^{2,\alpha}(\mathcal{M}\cup\mathcal N)$. In particular, \eqref{1388} holds in the classical sense. Finally,  the strong maximum principle and  the Hopf's lemma gives \eqref{142}.

\begin{flushright}
$\blacksquare$
\end{flushright}

The next result is used in Section \ref{210} to estimate in compact sets the solutions obtained above.

\begin{theorem} \label{14}  
Let  $X\subset \mathcal{M} \cup  \mathcal{N}$ be a compact set.  Suppose that $c, c_0 \in L^{\infty} (\mathcal{M})$,  $c_1, c_2 \in L^{\infty} (\mathcal{N})$ and $\|c\|_{L^\infty}  + \|c_1\|_{L^\infty}\leq \Lambda$. Let  $ 1<\alpha \leq \frac{n+2}{n-2} $  and  $ 1<\alpha _1 \leq \frac{n}{n-2}$.  Assume that  $u\in H^1 (\mathcal{M})$,  $u\geq 0$  and
\begin{equation} \label{16}
\int_{\mathcal{M}} \left(  g(\nabla u , \nabla v) + c u v +  c_0 u^{\alpha} v\right) \dv _g +\int_{\mathcal{N}} \left( c_2 u v  + c_1 u^{\alpha _1 } v \right) \s _g\leq 0, 
\end{equation} 
for all  $v\in H^{1}_{ \mathcal{D}}(\mathcal{M})$ with $v \geq 0$.
 
If  $c\geq -S$,  $c_2 \geq -S_2$, with  $S,  S_2 \geq 0 $,  and  $ c_0 \geq S_0$,  $c_1 \geq S_1$, with $S_0 ,  S_1 >0$,  then
\begin{equation} \label{17}
\stackbin[X]{}{\sup}  \ u \leq C\left(\Lambda , S , S_0 , S_1 , S_2 , \alpha , \alpha _1 , n , g ,X \right).
\end{equation}
\end{theorem} 
 
\noindent \textit{Proof.}  
We follow the steps of the proof of  \cite[ Theorem 1.1]{1aviles}. By the compactness of  $X$  we can find  $\varepsilon >0$  and a finite number of charts  $\left(  B_{3\varepsilon } , \varphi _i \right)$,  $\left(  K_{3\varepsilon} \cap \{x_n \geq 0 \} , \psi _j \right)$   such that  
$$
X \subset \left(\cup _i \varphi _i \left( B_\varepsilon  \right)\right) \cup \left( \cup _j \psi _j \left(  K_{\varepsilon} \cap \{x_n \geq 0 \} \right) \right),
$$
$\overline{\psi _j \left(  K_{3\varepsilon} \cap \{x_n = 0 \}  \right)} \subset  \mathcal{N}$ 
and  $\overline{ \varphi _i \left( B_{3\varepsilon} \right)} \subset \mathcal{M}$.  

We have the following two cases:
\begin{enumerate}[$(i)$]
\item $\stackbin[X]{}{\sup}  \ u = \stackbin[\varphi _i \left( B_\varepsilon \right) ]{}{\sup} \  u $,   for some  $i$.
\item $\stackbin[X]{}{\sup} \ u = \stackbin[\psi _j \left( K^+ _\varepsilon  \right) ]{}{\sup}  \ u $,  for some $j$,  where  $K^+ _\varepsilon   :=  K_\varepsilon \cap \left\{x_n > 0\right\}$.
\end{enumerate}
If the first case holds, the proof of   \eqref{17}  is the same as in \cite{1aviles}. For the second case, observe that, by \eqref{16},   
\begin{equation} \label{124}
\int_{\mathcal{M}} g(\nabla u , \nabla v) \dv_g -  S\int _{\mathcal{M}}  u v \dv _g   - S_2 \int _{ \mathcal{N}}  u v \s _g \leq 0,
\end{equation}
for all  $v\in H^{1}  _{\mathcal{D}} (\mathcal{M})$  with  $v\geq 0$. 

Set $K^0 _{\varepsilon}:=K_\varepsilon \cap \{ x_n =0 \}$ and let $p, p_1\geq 2$ be  constants to be determined below. Applying Lemma \ref{13} to inequality \eqref{124} we have
\begin{equation} \label{18}
\stackbin[K^+ _\varepsilon  ]{}{\sup} \ u\circ \psi _j + \stackbin[K^0 _\varepsilon  ]{}{\sup} \  u\circ \psi _j  \leq  C\left(S, S_2 , g , p, p_1 , n , K^+ _{3\varepsilon} \right) \left(  \Vert u\circ \psi _j )\Vert _{L^p \left(K^+ _{2\varepsilon}  \right)} + \Vert  u\circ \psi\Vert _{L^{p_1} (K^0 _{2\varepsilon}  )}  \right).  
\end{equation}

Now, let  $\zeta \in C_c ^\infty \left(\psi _j \left( K _{3\varepsilon}  \cap \{x_n \geq 0\}\right)\right)$, $\zeta \geq 0$, be such that $\zeta \equiv 1 $ on $\psi _j \left( K _{2\varepsilon}  \cap \{x_n \geq 0\}\right)$. Setting $\beta = \frac{2(\alpha+1) }{\alpha-1}>2$ and replacing $u\zeta ^\beta$ into \eqref{16}  we have
\begin{align*}
&S_0 \int_{\mathcal{M}} u^{\alpha +1} \zeta ^\beta \dv _g + S_1 \int_{\mathcal{N}} u^{\alpha _1 + 1} \zeta ^\beta \s _g \\ 
&\leq  -\int_{\mathcal{M}} \left[\zeta ^\beta |\nabla u|^2 + g\left(\nabla u,\beta u \zeta ^{\beta-1} \nabla \zeta \right) \right] \dv _g + \int_{\mathcal{M}} S u^{2} \zeta ^\beta  \dv _g  + \int _{\mathcal{N}} S _2 u^2 \zeta ^\beta \s _g .
\end{align*}
By Cauchy-Schwarz,
$$
-\beta u \zeta ^{\beta-1} g\left(\nabla u ,  \nabla \zeta \right) \leq \zeta ^{\beta}|\nabla u|^{2} + \beta^{2} u^{2} \zeta ^{\beta -2} |\nabla \zeta |^{2},
$$
so we obtain
\begin{equation*}
\begin{split}
&S_0 \int_{\mathcal{M}} u^{\alpha +1} \zeta ^\beta  \dv _g +  S_1 \int_{\mathcal{N}} u^{\alpha _1 + 1} \zeta ^\beta \s _g  
\leq   \beta^{2}  \int_{\mathcal{M}} u^{2} \zeta ^{\beta -2} |\nabla \zeta |^{2} \dv _g
+  S\int_{\mathcal{M}}  u^{2} \zeta ^\beta  \dv _g + S_2\int _{\mathcal{N}}  u^2 \zeta ^\beta \s _g,\\
& \leq   \beta ^2 \left( \int _{\mathcal{M}} u^{\alpha+1} \zeta ^\beta  \dv _g \right)^{\frac{2}{\alpha+1}}\left(\int _{\mathcal{M}} |\nabla \zeta |^\beta \dv _g \right)^{\frac{\alpha-1}{\alpha+1}} \\
&+ S \left(\int _{\mathcal{M}} u^{\alpha+1} \zeta ^\beta \dv _g \right)^{\frac{2}{\alpha+1}}  \left(\int _{\mathcal{M}}   \zeta ^\beta \dv _g \right)^{\frac{\alpha-1}{\alpha+1}} 
+ S_2 \left(\int _{\mathcal{N}} u^{\alpha _1 + 1} \zeta ^\beta \s _g \right)^{\frac{2}{\alpha _1 +1}}  \left(\int _{\mathcal{N}}   \zeta ^\beta  \s _g\right)^{\frac{\alpha _1 -1}{\alpha _1 +1}}.
\end{split}
\end{equation*}
Then, by Young's Inequality,
$$
 \int_{\mathcal{M}} u^{\alpha +1} \zeta ^\beta \dv _g + \int_{\mathcal{N}} u^{\alpha _1 + 1} \zeta ^\beta \s _g \leq C _1 \left(\int _{\mathcal{M}} |\nabla \zeta |^\beta \dv _g +  \int _{\mathcal{M}}   \zeta ^\beta \dv _g  + \int _{\mathcal{N}}   \zeta ^\beta \s _g \right),
$$
where  $C_1 =C_1 \left( S , S_0 , S_1 , S_2 , \alpha , \alpha _1 \right) >0 $. Then
$$
\| u \circ \psi _j \|_{L^{\alpha +1 } \left(K^+ _{2\varepsilon}\right)} + \| u \circ \psi _j \|_{L^{\alpha _1 + 1} \left(K^0 _{2\varepsilon}\right)} \leq C_2 \left(S, S_0 , S_1 , S_2 , \alpha , \alpha _1 , n ,  X \right). 
$$

Finally,  if   $p=\alpha+1$  and  $p_1 = \alpha _1 +1$, by \eqref{18} we have
$$
\sup _X \ u \leq C\left(\Lambda , S, S_0 , S_1 , S_2 , \alpha , \alpha _1 , n , g , X \right).
$$

\begin{flushright}
$\blacksquare$
\end{flushright}

\end{document}